\documentclass[12pt]{article}
\usepackage{amsfonts,amsmath,amsthm,amstext,latexsym,graphicx,bezier}
\usepackage{lineno}
\usepackage[matrix, arrow, curve]{xy}
\usepackage{enumerate}
\usepackage{subcaption}
\usepackage{epstopdf}
\usepackage[cp1250]{inputenc}
\usepackage{multirow}
\usepackage{booktabs}
\usepackage{colortbl}
\usepackage{tikz}
\usepackage{algorithm,algorithmic}
\numberwithin{equation}{section}
\usetikzlibrary{intersections,patterns,angles,quotes,calc}%%%%%

\newtheorem{propn}{Proposition}[section]

\newtheorem{theorem}{Theorem}[section]

\begin{document}

\begin{center}{\bf{\Large Sampling from Unknown Transition  Densities  of   Diffusion processes}}
\end{center}
%\begin{center}{\bf{\Large Rejection sampling  of paths  for   Diffusion processes via Fokker-Planck Equation}}
%\end{center}
Yasin Kikabi, Juma Kasozi\\%$^{1}$, \\
Department of Mathematics, College of  Natural Sciences, Makerere University.\\ P.O. Box: 7062 Kampala, Uganda.\\\\{\bf{Abstract}}\\
In this paper, we introduce a new  method of sampling from transition densities of diffusion processes including those unknown in closed forms by solving a partial differential equation satisfied by the quotient of transition densities. We demonstrate the performance of the developed method on processes  with known densities and the obtained results are  consistent with theoretical values. The method is applied to Wright-Fisher diffusions owing to their   importance in population genetics in studying   interaction networks  inherent in  genetic data. Diffusion processes   with bounded drift and non degenerate diffusion are considered as reference processes. \\ \\

$\bf {Key words}:$ Stochastic differential  equation (SDE), Transition density, Fokker-Planck  partial differential equation, Aronson's bound, Rejection sampling, Wright-Fisher diffusion.
\section{Introduction}

There are several  ways  of representing stochastic systems mathematically  for instance   using stochastic differential equations.  A stochastic differential  equation (SDE) is an equation of the form 
\begin{eqnarray}
\label{eqn1}
dX_t=a(X_t,t)dt+b(X_t,t)dW_t,
\end{eqnarray} where  $a(X_t,t)$ is the drift term, $b(X_t,t)$ is diffusion coefficient and $W_t$ is the Wiener process.
Due  to the interplay  between stochastic differential   equations and the Fokker- Planck equation see (Gardiner 1985, chapter $5$)\cite{gardiner1985handbook}, the transition density of the process described by (\ref{eqn1}) is obtained by solving the Fokker-Planck equation associated with the SDE.  
The Fokker-Planck equation is a partial differential equation that describes the evolution of the conditional probability distribution say $ P(x, t/x_0,t_0)$ of a stochastic system  through time in terms of state-space variables   $x$. 
For instance, following Grigorios (2014)\cite{pavliotis2014stochastic},  the one dimensional,  the Fokker-Planck equation associated with (\ref{eqn1}) is given by 
\begin{eqnarray}
\label{eqn2}
\frac{\partial}{\partial{t}}{P(x,t)}=-\frac{\partial}{\partial{x}}[a(x,t)P(x,t)]+\frac{1}{2}\frac{\partial^2}{\partial{x^2}}{[b(x,t)P(x,t)]}
\end{eqnarray}  with initial condition $ P(x, t_0/x_0,t_0)=\delta(x-x_0)$.

In equation (\ref{eqn2}),  $P(x,t)$ denotes the conditional density  $ P(x, t/x_0,t_0)$. 

Unfortunately, obtaining closed form solutions to   Fokker-Planck equations associated with many  diffusion processes of practical interest  such as the Wright-Fisher  diffusion  remains a challenging problem Yun \& Steinrücken (2012 )\cite{song2012simple}. This has always hampered  the likelihood based inference for  diffusion models with discretely observed data due to the fact   it requires  the transition density to be known explicitly which is hardly the case.  However,  several methods have been proposed in literature to  circumvent the intractable likelihoods in statistical inference of diffusion models described by SDEs. These  include  approximating the transition density $P$ by a Gaussian density (Florens-Zmirou  1989; Yoshida 1992; Sorensen \& Uchida 2003)\cite{florens1989approximate, yoshida1992estimation, sorensen2003small},  approximating $P$ by orthogonal polynomials, Ait-Sahalia  (2002, 2008)\cite{ait2002maximum, ait2008closed} and approximating  $P$ by solving the Fokker-Planck Equation numerically Bollback et al. (2008)\cite{bollback2008estimation}.\\ On the other hand, alternative  inference methods such as  the so-called Approximately Bayesian Computation(ABC)  for diffusion models found for instance in Tina at el. (2009)\cite{toni2009approximate}, require    obtaining  exact sample paths from the diffusion models. That is,  sampling exactly  from the transition densities associated with the SDEs describing the process. However, these transition densities  are mostly  known by eigen-fuction expansion or  by their   Fokker-Planck equations  which itself is  a challenge to solve to obtain   closed forms of transition densities.

Recently,  exact simulation methods for diffusion process have been proposed by  Beskos et al.(2005, 2006, 2008) \cite{beskos2005exact,beskos2006retrospective, beskos2008factorisation}. These are based on equivalence of   measures induced by  SDEs that comprise the   same diffusion coefficient where rejection sampling is obtained via the Radon-Nikodym derivative. These methods  require  obtaining  exact sample paths from one of the diffusions which are used as proposal paths in the rejection sampling.  Generally, Beskos's utilizes the absolute continuity of measures induced by diffusions  with the same diffusion coefficient  via Girsanov's theorem to carry out rejection sampling of diffusion paths.\\
In this work, we explore the possibility of obtaining sample paths for diffusions by  rejection sampling  via Fokker -Planck equations.  We consider diffusions that have the same diffusion coefficient.  Our method is based on the assumption that the transition density of one of the diffusions (proposal diffusion) is known so that we can obtain samples from it  at any discrete time points. This method involves solving a partial differential equation that involves quotients of transition densities of the target diffusion and proposal diffusion. Each of the  two diffusion processes  is   defined  by a system of stochastic differential equations.\\ Differently from Beskos's method, where the entire path is either accepted or rejected, a path is developed by sampling from transition densities at discrete  time points using the classical rejection sampling from probability density functions. This method is based on the assumption that we can solve the PDE that  involves quotients of  conditional transition densities.

The remainder of this paper is organised as follows. In Section $2$ we present the classical rejection sampling method for purposes of development of our method.  In Section $3$, we present our rejection method, and its implementation.   In Section $4$  we apply our method to  diffusion processes with bounded drift and non degenerate diffusion coefficients that satisfy Aronson's bound with specific application to the  Wright-Fisher diffusion. Later, we consider processes with unbounded drift . We discuss conclusions and possible directions for future research in the last section.
  \section{General rejection sampling} 
  Rejection sampling is a widely applied technique of sampling from one probability density using samples from another density.  Following Beskos ( 2005)\cite{beskos2005exact}, the classical rejection method is presented as follows.\\
Assume that $P_1$ and $P_2$ are probability densities with respect to  some measure on $\mathbb{R}^d $ $d>0$ and that there exists $\epsilon>0$ such that $\frac{P_2}{\epsilon P_1} \leq1$.Then, the following algorithm returns samples distributed according to  $P_2$.
 \begin{algorithm}
\caption{Rejection   algorithm }
\label{algo_1}
\begin{algorithmic}[1]
\STATE SAMPLE Y$\sim$$P_1$.

\STATE  SAMPLE U$\sim $ Unif(0,1).
 \STATE IF U$<\frac{P_2}{\epsilon  P_1}(Y)$  
\RETURN Y 
\STATE ELSE go to $1$

\end{algorithmic}
\end{algorithm}    
Observe that, the rejection condition requires obtaining values of the quotient of probability densities at the sample points $Y$. Obtaining such values is challenging specifically in cases where the closed form of  $P_2$  is unknown.  In our method presented in the next section, obtaining such values involves solving a partial differential equation of quotients of probability densities so that (Algorithm \ref{algo_1} ) is then applicable.

     \section{Rejection sampling via the Fokker-Planck equation } 
 Consider the following multidimensional  stochastic differential equations    
\begin{eqnarray}
\label{eqn3}
d{\bf{X}}_{t}&=&{\bf{S}}_i({\bf{X}}_t)dt+\sqrt{\sigma({\bf{X}}_t)}d{\bf{W}}_t, \,\,\, {\bf{X}}_0,\,\, 0\leq t \leq T,\,\,\, i=1,2,
\end {eqnarray} 
where  ${\bf{S}}_i$ , ${\bf{X}}_t\in \mathbb{R}^d$, $d\in\mathbb{N}$ and ${\sigma({\bf{X}}_t)}$ is a $d-$square matrix.

Let  ${P}_1$ and ${P}_2$ be  the corresponding conditional probability  densities with $P_1$ known explicitly in its closed form.  Then the corresponding Fokker-Planck equations are given by 
\begin{eqnarray}
\label{eqn41}
\frac{\partial}{\partial{t}}{P_i({\bf{x}},t)}=-\nabla_{\bf{x}}^T({\bf{S}}_i({\bf{x}})P_i)+\frac{1}{2}\nabla_{\bf{x}}^T(\sigma({\bf{x}}) P_i)\nabla_{\bf{x}},\,\,\, 
\end{eqnarray} 
  subject to initial condition $P_i({\bf{x}},t_0)=\delta({\bf{x}}-{\bf{x}}_0)\,\,\,i=1,2.$ 
   $\nabla_{\bf{x}}^T$ denotes  the multivariate differential operator, 
 $\nabla_{\bf{x}}^T=(\frac{\partial}{\partial{x_1}}, \frac{\partial}{\partial{x_2}},...,\frac{\partial}{\partial{x_n}})$.

 {\bf{Assumptions:}} 
 Let us introduce the following set of assumptions, 
\begin{itemize} 
\item[{\bf{A1}.}] the Fokker-Planck equation for $P_2$ yields a  solution on time interval  $[0,T]$ such   there exists $c\in \mathbb{R}$, $c>0$ for which  $\frac{P_2({\bf{x}},t)}{P_1({\bf{x}},t)}\leq c$ for any ${\bf{x}}\in \mathbb{R}^d$ and $t\in[0,T]$. 
 \item[{\bf{A2}.}]  ${\bf{S}}_1({\bf{x}})$, ${\bf{S}}_2({\bf{x}})$ and $\sigma_{k,j}({\bf{x}}),\,\,\, j,k=1,2,...,d.$  are differentiable functions. 
  \end{itemize} Then, the quotient of transition densities $\frac{P_2}{P_1}({\bf{x}},t)$ satisfies a  partial differential equation given in the following theorem.

 { \theorem
  Let \begin{eqnarray*}
\label{eqn31}
d{\bf{X}}_{t}&=&{\bf{S}}_i({\bf{X}}_t)dt+\sqrt{\sigma({\bf{X}}_t)}d{\bf{W}}_t, \,\,\, {\bf{X}}_0,\,\, 0\leq t \leq T,\,\,\, i=1,2.
\end {eqnarray*} be SDEs satisfying assumptions (A$1$) and (A$2$) above. Then the quotient of transition densities  $\frac{P_2}{P_1}({\bf{x}},t)$ satisfies  the partial differential equation

 \begin{eqnarray} 
   \frac{\partial}{\partial{t}}{\big[\frac{P_2}{P_1}({\bf{x}},t)]=\big\{\nabla_{\bf{x}}^T({\bf{S}}_1-{\bf{S}}_2)+\nabla_{\bf{x}}^T \log(P_1)({\bf{S}}_1-{\bf{S}}_2)\big \}[\frac{P_2}{P_1}] }+\nonumber \\ \big\{(\nabla_{\bf{x}}^T \log(P_1))\sigma({\bf{x}}) } +\nabla_{\bf{x}}^T\sigma({\bf{x}})- {{\bf{S}}_2    \big\}^T\nabla_{\bf{x}}[\frac{P_2}{P_1}] +\frac{1}{2}trace
   \{ \sigma({\bf{x}})\nabla_{\bf{x}}(\nabla_{\bf{x}}[\frac{P_2}{P_1}])\}. 
   \end{eqnarray}
   $A^T$ denotes the transpose of the matrix $A$.

\begin{proof}
We sketch the proof for $d=1$, the general case being straight forward.
Let d=1 and denote  the second partial derivative of $P$ with respect to $x$ by $\partial^2_{xx}P$. 

Consider
\begin{eqnarray}
\frac{\partial}{\partial{t}}{\big[\frac{P_2}{P_1}(x,t)\big]}=\frac{P_1{\partial_tP_2}-P_2 {\partial_tP_1}}{P_1^2}.
\end{eqnarray}
Using the Fokker-Planck  equations associated with the constituent SDEs  equation (\ref{eqn41}),  yields
\begin{eqnarray}
\label{44}
\frac{\partial}{\partial{t}}{\big[\frac{P_2}{P_1}(x,t)\big]}=\frac{-P_1{\partial_x(S_2(x)P_2)}+P_2{\partial_x(S_1(x)P_1)}}{P_1^2}+\frac{P_1\partial_{xx}^2(\sigma(x)P_2)-{P_2}\partial_{xx}^2(\sigma(x)P_1)}{2P_1^2}.
\end{eqnarray}
Note that the second order partial derivatives in equation (\ref{44}) can be expressed as  
\begin{eqnarray}
\label{424}
\partial_{xx}^2(\sigma(x)P_2)=\sigma(x)\partial_{xx}^2P_2+2\partial_{x}\sigma(x)\partial_xP_2+\partial_{xx}^2\sigma(x)P_2,\end{eqnarray} and similarly, 
\begin{eqnarray}
\label{434}
\partial_{xx}^2(\sigma(x)P_1)=\sigma(x)\partial_{xx}^2P_1+2\partial_{x}\sigma(x)\partial_xP_1+\partial_{xx}^2\sigma(x)P_1.
\end{eqnarray}
By equation (\ref{424}) and (\ref{434}), 
\begin{eqnarray}
P_1\partial_{xx}^2(\sigma(x)P_2)-{P_2}\partial_{xx}^2(\sigma(x)P_1)=\sigma(x)(P_1\partial_{xx}^2P_2-P_2\partial_{xx}^2P_1)+2\partial_{x}\sigma(x)(P_1\partial_xP_2-P_2\partial_xP_1).
\end{eqnarray}
Thus, the last term in equation (\ref{44}) simplifies to 
\begin{eqnarray}
\frac{\sigma(x)}{2P_1^2}\partial_x(P_1^2\partial_x\big[\frac{P_2}{P_1}\big])+\partial_{x}\sigma(x)\partial_x\big[\frac{P_2}{P_1}\big], 
\end{eqnarray}
which further reduces to 
\begin{eqnarray}
\label{45}
 \frac{1}{2}\sigma(x) \partial^2_{xx}\big[\frac{P_2}{P_1}\big] +(\sigma(x)\frac{\partial_xP_1}{P_1} +\partial_{x}\sigma(x))\partial_x\big[\frac{P_2}{P_1}\big].
 \end{eqnarray}
 On the other hand,  the remaining term on the right side of equation (\ref{44})  reduces to 
 \begin{eqnarray}
 \label{46}
-S_2(x)\{\partial_x\big[\frac{P_2}{P_1}\big]+\frac{P_2}{P_1}\frac{\partial_xP_1}{P_1}\} +S_1(x)\frac{P_2}{P_1}\frac{\partial_xP_1}{P_1} -(\partial_{x}S_2(x)-\partial_{x}S_1(x))\frac{P_2}{P_1}. \end{eqnarray}
Using  (\ref{45}) and (\ref{46}) yields the partial differential equation( \ref{eqn6}) which is satisfied by the ratio of transition densities,
\begin{eqnarray}
\label{eqn6}
\frac{\partial}{\partial{t}}{V(x,t)}=\Big\{ \partial_{x}[S_1(x)-S_2(x)]+[S_1(x)-S_2(x)]\partial_x\log {P_1}\Big \}V(x,t)+\nonumber \\ \Big\{\sigma(x)\partial_x\log {P_1}+\partial_{x}\sigma(x)-S_2(x) \Big\}\frac{\partial}{\partial{x}}{V(x,t)} +\frac{1}{2}\sigma(x)\frac{\partial^2}{\partial{x^2}}{V(x,t)},\end{eqnarray}  with initial condition $ V(x,t_0/x_0,t_0)=1$,  where  $V(x,t)=\frac{P_2}{P_1}(x,t)$. 
\end{proof}}

Our rejection sampling method is based on  Theorem $3.1$. To enhance  our method,  it requires establishing  existence of solution to the partial differential (\ref{eqn6}) which we briefly consider in Subsection $3.1$.   In application  of theorem $3.1$, we restrict our consideration  to a  class of multidimensional diffusion processes that satisfy  Aronson's bound and later to diffusions with unbounded drift term  in the one dimensional case. Before that,  we state  the following theorem due to Aronson (1967)\cite{aronson1967bounds}. This theorem will be vital in ensuring that assumption $\pmb{A}_1$ to Theorem $3.1$ holds in our application.
{\theorem{Aronson's bound}\\
Consider the stochastic differential equation
 \begin{eqnarray}
\label{eqn32}
d{\bf{X}}_{t}&=&{\bf{S}}({\bf{X}}_t)dt+{\sigma({\bf{X}}_t)}d{\bf{W}}_t, \,\,\, {\bf{X}}_0,\,\, 0\leq t \leq T,\,\,\, 
\end {eqnarray}   with state space $S$.  Assume there exists a constant $\lambda>0$ such that $\lambda ^{-1}I\leq {\sigma({\bf{x}})} {\sigma({\bf{x}})}^T\leq\lambda I $ 
for all  ${\bf{x}}\in S$ and that ${\bf{S}}({\bf{x}})$ is bounded over $S$. Then the bound to the transition density $p(x,s;y,t)$ associated with the SDE is given by 
$$K_1g(x,s;y,t)\leq p(x,s;y,t)\leq K_2g(x,s;y,t),$$ where  $K_1$, $K_2$ are constants and $g$ is the transition density  of Brownian motion}.
\subsection{Numerical solution to the partial differential equation satisfied by ratio of transition densities}  

By Theorem $3.1$, the partial differential equation (\ref{eqn6}) satisfied by ratio of transition densities is 
 is of a general form 
\begin{eqnarray}
\label{eqn17}
\frac{\partial}{\partial{t}}{V(x,t)}=c(x,t)V(x,t)+b(x,t)\frac{\partial}{\partial{x}}{V(x,t)} +a(x,t)\frac{\partial^2}{\partial{x^2}}{V(x,t)},
\end{eqnarray}  
with initial condition $ V(x,t_0)=1$. This is a second order, variable coefficient, linear parabolic partial differential equation  whose  approximate solutions  are  obtained by numerical methods such as finite difference methods.  Obtaining numerical solutions to (\ref{eqn17}) over the time interval $[0, T]$ and spatial interval $[x_{min}, x_{max}]$ leads  to   solving   an initial value boundary value PDE numerically.   $x_{min}$ and  $ x_{max}$ denote the  minimum and maximum sample values along the path obtained from $P_1$ respectively. 
\\In this work, we apply  the Cranck-Nicolson finite difference scheme to solve the developed PDE  of   quotients of densities. This is largely attributed to   its unconditional stability and consistency properties Marques (2017)\cite{marques2017option}. In the multidimensional case (Section $4.2$), the Alternating Direction Implicit (ADI) method is used. \\ Following Marques (2017)\cite{marques2017option}, for  a general PDE of the form (\ref{eqn17}), consider  mesh points $(x_j,t_i)$ with spatial and time spacing $h$ and $k$ respectively. That is, $x_j=jh$, $t_i=ik$ $i=1,2,...,N$, $j=1,2,..., M$; $M,N\in \mathbb{N}$. k Let  $u_j^i =u(x_j,t_i)$, $a_j^i = a(x_j,t_i)$, $b_j^i = b(x_j,t_i)$, $c_j^i = c(x_j,t_i)$.  The Cranck-Nicolson scheme is given by 
\begin{eqnarray}
\label{eqn18}
A_iu_{j+1}^{i+1} + B_iu_{j}^{i+1} + C_iu_{j-1}^{i+1} = D_i ,\end{eqnarray}
where $A_i = 2ka_{j}^{i+1} +hkb_j^{i+1}$, $ B_i = - 4h^2 - 4ka_j^{i+1} +2h^2kc_j^{i+1}$, $C_i = 2ka_j^{i+1 }- hkb_j^{i+1}$, \\$Di =  -[u_{j+1}^i-(2ka_j^i +hkb _j^i)+u_j^i(4h^2-  4ka_j^i +2h^2c_j^i)+u^i_{j-1}(2ka_j^i -hkb_j^i)]$, $i=1,2,...,N$, $j=1,2,..., M$. \\
$N$ and $M$ are  the number of subintervals time in time  and space respectively.\\
To obtain the solution at the next time line, we  solve  linear system given by (\ref{eqn18}).

Having  described how to obtain  values of $\frac{P_2}{P_1}(x_i,t_i)$ via the developed PDE,  we now describe how to obtain  sample path points from $P_2$ using exact sample paths from $P_1$.
The idea of obtaining such samples from $P_2$ is as  follows:\\ 
Obtain  sample path points from $P_1$ at time instances $(t_1,t_2,...,t_T)$. Denote the minimum observation by   $x_{min}$ and similarly the maximum observation by $x_{max}$.\\
Solve the  PDE (\ref{eqn6}) numerically at time instances $(t_1,t_2,...,t_T)$ subject to boundary conditions $V(x_{min},t)=1=V(x_{max},t)$. The $x$-mesh is strategically chosen such that  it contains both $x_{min}$ and $x_{max}$  and  all the   sample points.   This solution gives the approximate value of $\frac{P_2}{P_1}$ at sample points and thus the usual rejection sampling procedure of using uniformly generated random numbers on unit interval.  That is,  a sample point is  rejected, if  random number generated uniformly on the unit interval falls above $\frac{P_2}{P_1}(x_i,t_i)$, otherwise it is accepted as a sample point from $P_2$. This  rejection step can be carried out simultaneously  for all sample path points  by generating $T$ uniform random numbers. The values at any of the rejected points can be obtained by interpolation between any two accepted points using a diffusion bridge. The  sampling procedure is summarised in the following algorithm.

 \begin{algorithm}
\caption{Rejection sampling  algorithm }
\label{find_cltree}
\begin{algorithmic}[1]
\STATE  Sample paths points $x_i$ from $P_1$  at time times  $t_i$  $i=1,2,...$,  $ T\in \mathbb{N}$ 

\STATE Approximate the solution to the  PDE (\ref{eqn6}) at time instances $\{t_1,t_2,...,t_T \}$  with x-grid that includes all the sample points with boundaries at $ x_{min}$ and $x_{max}$
 \RETURN $\frac{P_2}{P_1}(x_i,t_i)$  $i=1,2,...,T$
\STATE Generate $T$ uniformly distributed numbers on $[0,1]$

\STATE For $i=1,2,...T$, if $u_i\leq \frac{P_2}{c P_1}(x_i,t_i)$, accept $(x_i,t_i)$  else reject  $(x_i,t_i)$

\RETURN  Accepted sample path points $(x_i,t_i)$. 
\end{algorithmic}
\end{algorithm}    
 \section {Numerical Results}
 \subsection{Diffusions with bounded Drift:}
\subsubsection{The $1$-dimensional Wright-Fisher diffusion   with selection.}
The Wright-Fisher diffusion is a continuous time Markov process with a continuous state space that is used to model fluctuations in size of finite populations under the influence of evolutionary forces such as selection and  mutations. In population genetics, the Wright-Fisher diffusion in form of a stochastic differential equation is used to model the changes in allele frequencies in finite populations. For instance, considering   a population under selection evolutionary force,  the Wright-Fisher diffusion in the univariate case is given by 
\begin{eqnarray}
       \label{eqn350}
   dx_t=\gamma x_t(1-x_t)dt+\sqrt{x_t(1-x_t) } dW_t 
        \end{eqnarray} 
where $x_t$ is allele frequency at time $t$ and $\gamma$ is the selection parameter.  Observe that the drift function is bounded on the state space $[0,1]$.
Data on such diffusion processes is increasing becoming available  due to advancement in biotechnology as earlier mentioned in the introduction of chapter $1$ with a challenge being the development reliable statistical methods that can accurately determine the parameters underlying such processes such as selection  and  mutation parameters in the univariate cases and interaction parameters in the multivariate cases.

Such statistical methods in cases of (\ref{eqn350}) has always been hampered by lack of closed form transition densities associated with the diffusion to explicitly express the likelihood function. However, lack of a likelihood function can be circumvented by using likelihood free inferential methods such the Approximately Bayesian methods.  
These rely on the possibility of simulating exactly from models which in this cases implies the need for simulating sample paths of the diffusion (\ref{eqn350}).
In this section, we show how to sample from the transition density of (\ref{eqn350})  using the transition density of a derived process from a Brownian bridge  using  our method.

Consider   Brownian motion   $x_t$ on the time interval $[0,T]$  started at $x_0$. This  is described  by the  stochastic differential equation, 
 \begin{eqnarray}
 dx_t=dW_t,\,\,\, x(0)=x_0. 
 \end{eqnarray}

 Take the transformation of $x_t$ under the map
\begin{eqnarray}
\label{sigmoid}
     g(x)=\frac{1}{1+e^{-x}},
        \end{eqnarray} 
      and denote the derived process by $y_t$.
       By Ito's Lemma, the derived  process $y_t$  satisfies the following  stochastic differential equation
     \begin{eqnarray}
     \label{eqn_21}
            dy_t=\frac{1}{2}y_t(1-y_t)(1-2y_t)dt+{y_t(1-y_t) } dW_t, \,\,\, y_0=Y_0,\,\,\, 0<t<T,
        \end{eqnarray} 
where  $W_t$ is the standard Wiener process.
       The derived process has a state space  $(0,1)$  as the Wright-Fisher diffusion assuming inaccessible boundaries  at $0$ and $1$. This  suggests  the derived process $y_t$ as  a good candidate for the reference process in our rejection sampling for  the Wright-Fisher diffusion if its transition density is known in a closed form.\\
       However, noting that $y_t$ is a  bijective  transformation of a Gaussian process under $g$, its transition density can be obtained by the transformation Theorem which we state shortly  in the following paragraph,  a detailed proof of which is found   for example in Gut (  2009 )\cite{gut2009intermediate} Chapter $1$. Following Gut (2009 )\cite{gut2009intermediate}, the transformation theorem is stated as follows, 
      {\theorem {(Transformation theorem)}\\
      Let ${\bf{X}}$ be an $n$-dimensional, continuous, random variable with density $f_X(x)$, and suppose that ${\bf{X}}$  has its mass concentrated on a set $S\subset \mathbb{R}^n$. Let $ g = (g_1, g_2, ..., g_n)$ be a bijection from $S$ to some set $T\subset \mathbb{R}^n$, and consider the $n$-dimensional random variable  $${\bf{Y}}=g({\bf{X}}).$$  
      Assume that g and its inverse are both continuously differentiable.\\
      The density of  ${\bf{Y}}$ is given by 
 \begin{eqnarray} 
          \label{eqn_33}    
        f_{\bf{Y}}({\bf{y}})=
        \begin{cases}
      & f_{\bf{X}}(h_1({\bf{y}}),h_2({\bf{y}}),\cdots,h_n({\bf{y}}))|{{\bf{J}}}| \,\,\, for\,\,\, {{\bf{y}}\in T}\\
      &0 \,\,\,\,\,\,\,\,\text{otherwise},
        \end{cases}        
\end{eqnarray}      
where $h$ is the inverse of $g$ and $|J|$ is the Jacobian.}
\\
Denote the transition density of (\ref{eqn_21}) by $P_1$. By the transformation theorem $4.1$,  the transition density of  the derived  process  equation (\ref{eqn_21}) is obtained from the Gaussian density of a Brownian bridge  and thus,  given by 
  \begin{eqnarray} 
    P_1(s,y_s, t,y_t)=\frac{1}{\sqrt{2\pi (t-s)}}\frac{e^{-\frac{1}{2}\frac{(\beta_t-\beta_s)^2}{t-s}}}{y_t(1-y_t)},\,\,\, y_t\in(0,1),
    \end{eqnarray} 
    where $\beta_t=\log{(\frac{y_t}{1-y_t})}$.
    
    We now embark on sampling from the Wright-Fisher process using the process  (\ref{eqn_21}) as a reference process.      Note that   both processes share the same state space, however, their  diffusion terms are different. \\      
   Thus, to  employ our method, the diffusion (\ref{eqn350}) has to be transformed into a diffusion with the same diffusion coefficient  as equation (\ref{eqn_21}). In this regard , consider  the map
       \begin{eqnarray} 
          \label{eqn_34}    
        f(x)=\frac{1}{1+e^{-\sin^{-1}(2x-1)}}.
        \end{eqnarray}       
 Under this transformation, (\ref{eqn350}) is transformed into an Ito process given by the following stochastic differential equation,
  \begin{eqnarray}
     \label{eqn_4}
     dy_t=\frac{1}{2}y_t(1-y_t)(1-2y_t-\tan \beta_t +\gamma \cos(\beta_t))dt+{y_t(1-y_t) } dW_t, \,\,\,y_0=Y_0     
    \end{eqnarray}   where $\beta_t=\log(\frac{y_t}{1-y_t})$.

 Now, we use the  process given by  equation (\ref{eqn_21}) as the reference process in our sampling algorithm to obtain  samples from the transition density of the stochastic  process given by  equation (\ref{eqn_4}). The obtained samples are later transformed into samples from the Wright-Fisher process  equation (\ref{eqn350}) by the inverse map $f^{-1}(x)$.
 
 Observe that both processes (\ref{eqn_4}) and  (\ref{eqn_21}) have bounded drift functions on (0,1) and the diffusion coefficient is non degenerate over the same set. Hence, the two processes satisfy Aronson's bound (theorem $3.2$) so that the  assumption ${\bf{A}}1$ holds.
 
 Thus, by  theorem $3.1$, the  partial differential equation satisfied by the ratio of transition densities is given by 
 
 \begin{eqnarray}
\label{eqn_5}
\frac{\partial V}{\partial{t}}&=&\frac{1}{2}(\gamma \sin\beta+\sec^2\beta+\frac{(\beta-\beta_0)(\gamma \cos\beta-\tan\beta)}{t})V
    +\nonumber \\ &&\frac{1}{2}y(1-y)(1-2y+2\frac{\beta_0-\beta}{t}+\tan\beta-\gamma\cos\beta)\frac{\partial V}{\partial{y}} +
\nonumber \\&&\frac{(y(1-y))^2}{2}\frac{\partial^2 V}{\partial{y^2}},\,\,\, t>0,
\end{eqnarray}  
where  $\beta(y)=\log{(\frac{y}{1-y})}$, $V(y,t)=\frac{P_2}{P_1}(y,t)$ and $P_2$ is the transition density to the transformed Wright-fisher Process (\ref{eqn_4}) . The  PDE is subjected to initial conditions $\frac{P_2}{P_1}(y,0)=1$  and boundary conditions  $\frac{P_2}{P_1}(y_{min},t)=1$ ,  $\frac{P_2}{P_1}(y_{max},t)=1.$  \\ The simulated solution to the PDE is shown in Figure  $1$ with a corresponding accepted path in Figure $2$.
\begin{figure}[h]
  \centering
  \includegraphics[ width=0.4\textwidth]{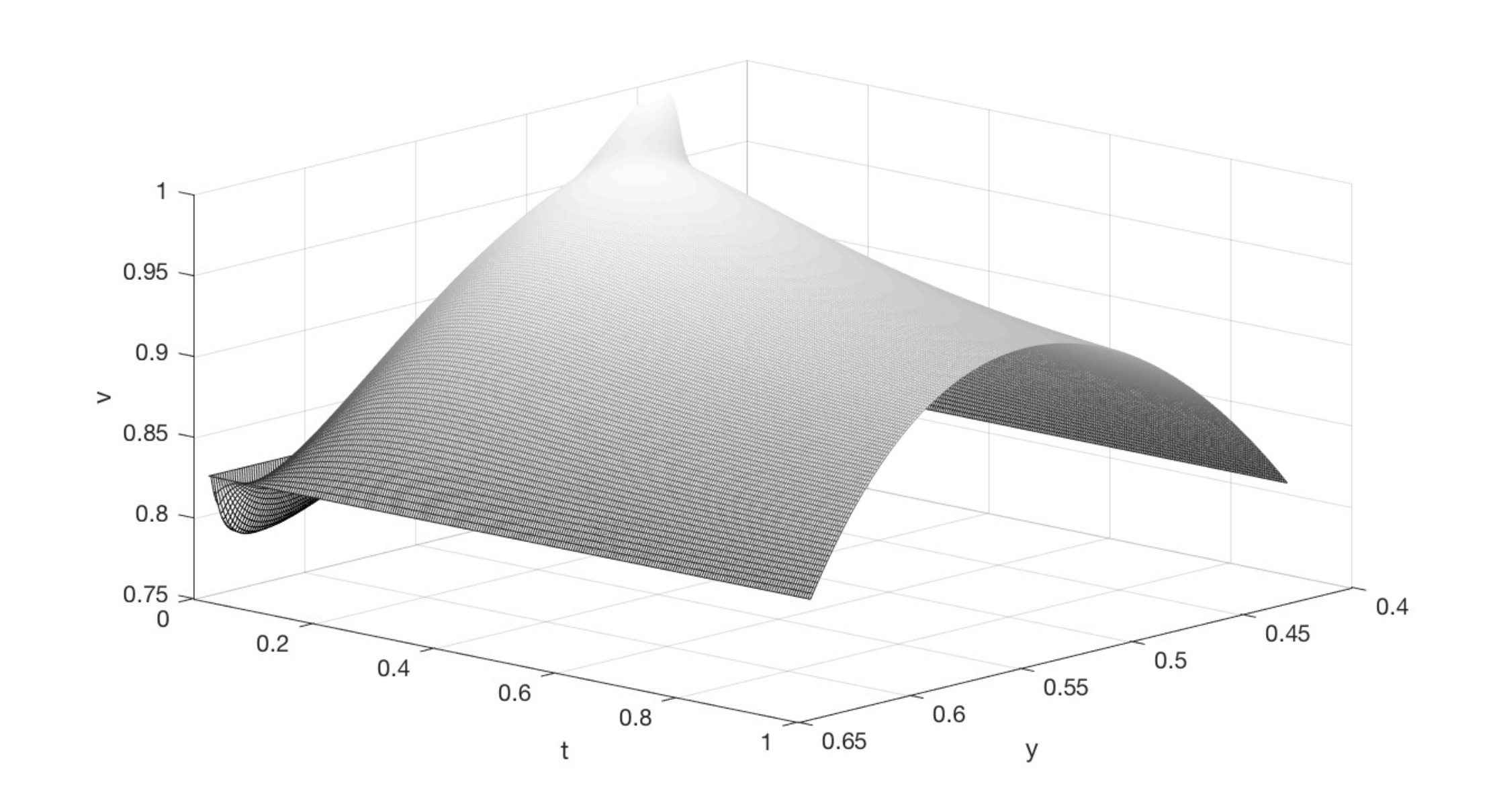}

 \caption{{\bf{Fig.1}} Numerical solutions of the PDE for $\gamma=1$ at points on the proposed transformed  Brownian path }
  \label{fig:31}
\end{figure}
  \begin{figure}[h]
  \centering
  \includegraphics[ width=0.43\textwidth]{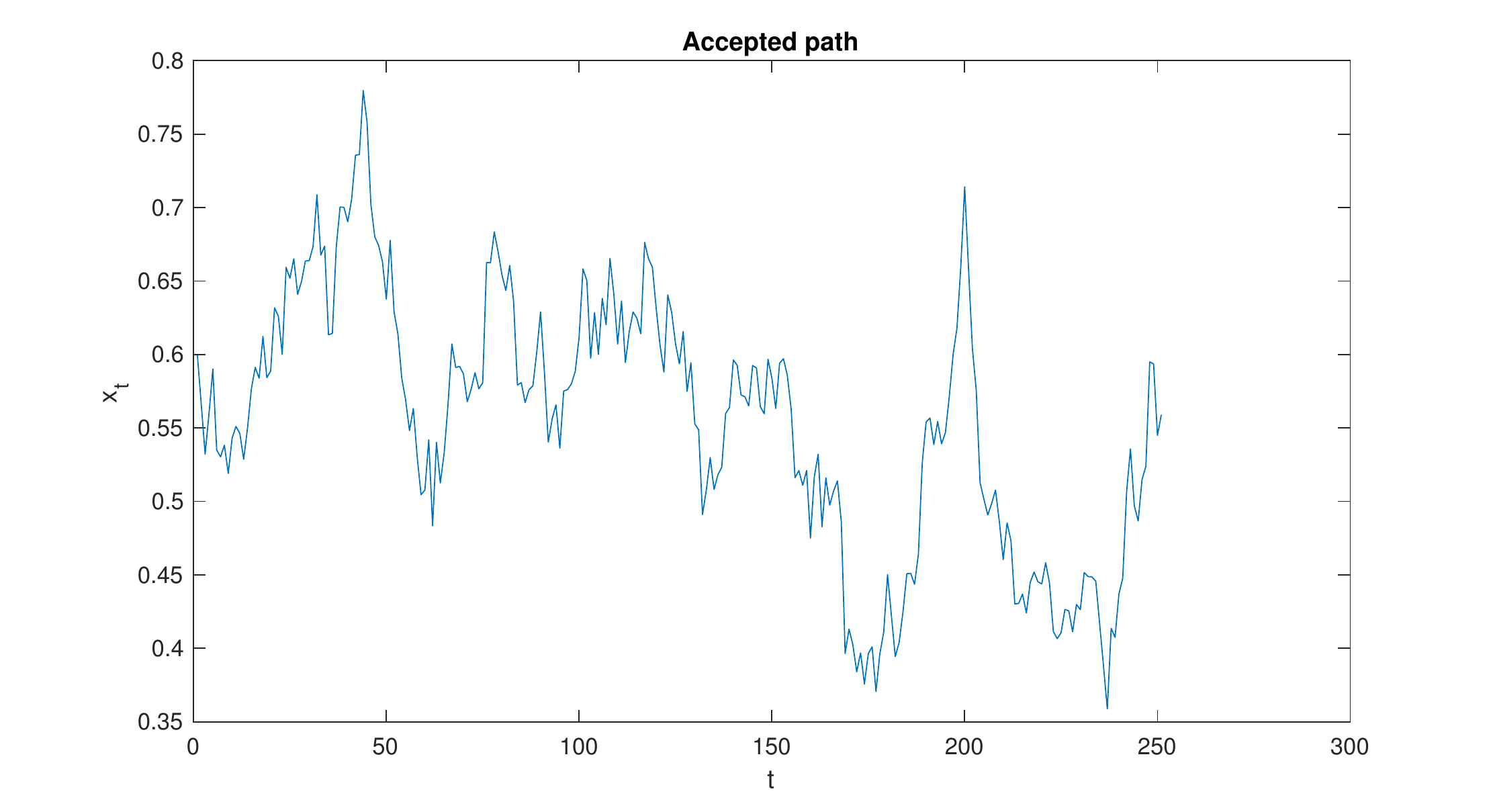}
    \caption{{\bf{Fig.2}} Accepted  transformed Brownian motion path  }
  \label{fig:32}
\end{figure}
\newpage
In addition, we observe that paths that spend a significant time near the right  boundary are  rejected   since  PDE yields mostly zero values including  negative values at some points  that carry no physical meaning,  see for instance Figure $3$.   
  \begin{figure}[h]
  \centering
  \includegraphics[ width=.43\textwidth]{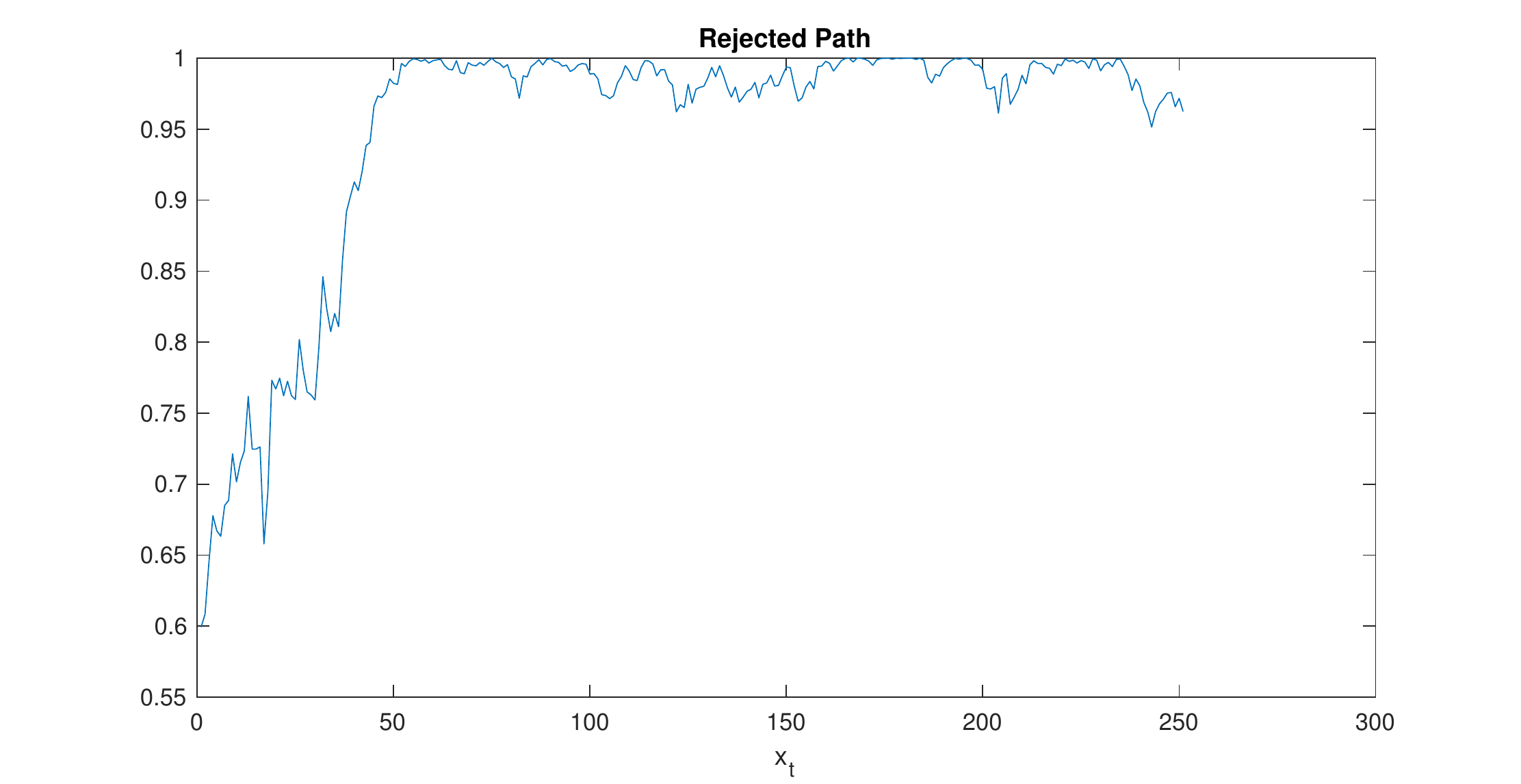}
  \includegraphics[ width=.43\textwidth]{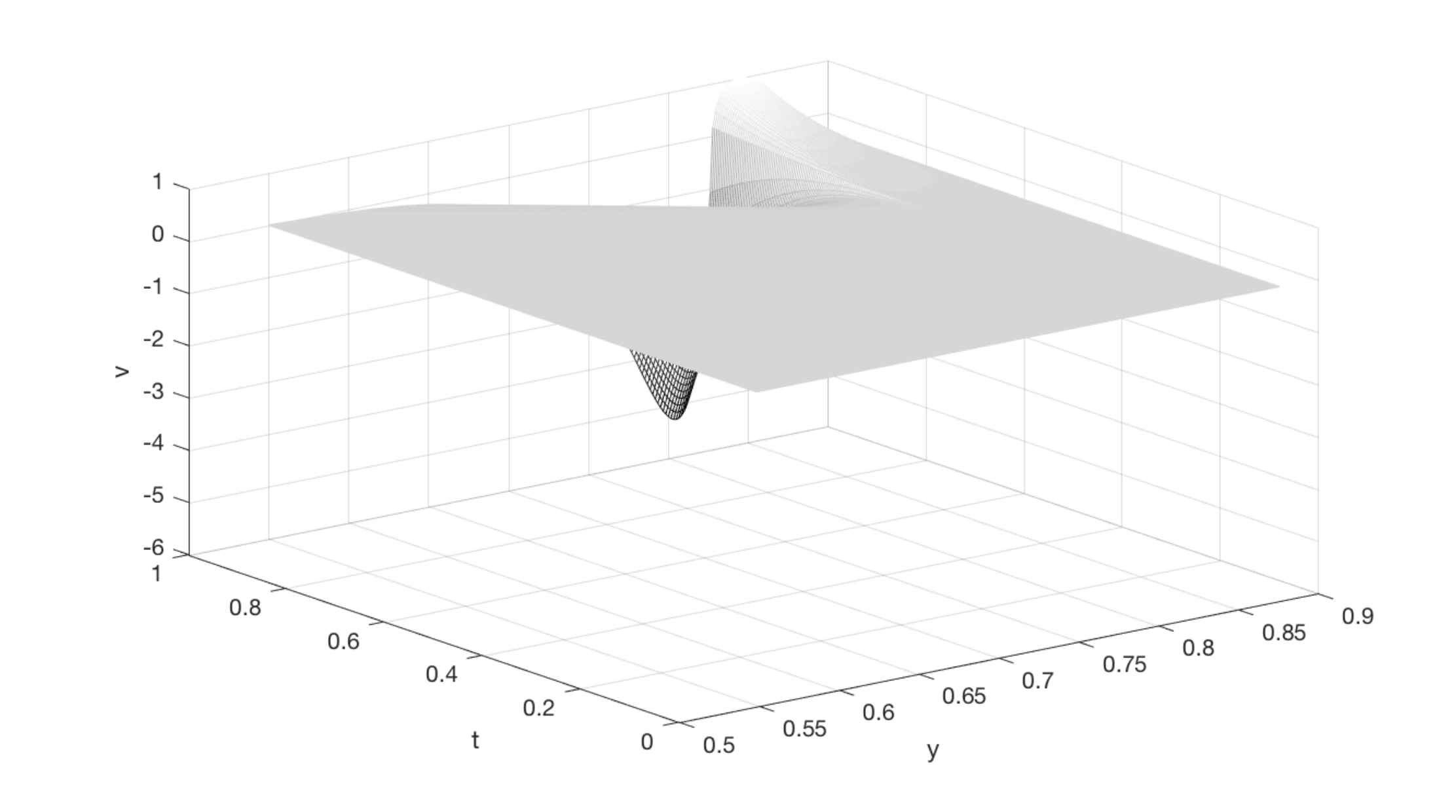}
  \caption{{\bf{Fig.3}} Rejected path with corresponding numerical solutions of the PDE for $\gamma=1$ }
  \label{fig:351}
\end{figure}  

Since statistical inference of diffusion processes  assumes observed  data, the accepted paths  in our case are transformed into Brownian bridges. Figure $4$ shows the Brownian bridges corresponding to accepted paths.
  \begin{figure}[h]
  \centering
  \includegraphics[ width=0.43\textwidth]{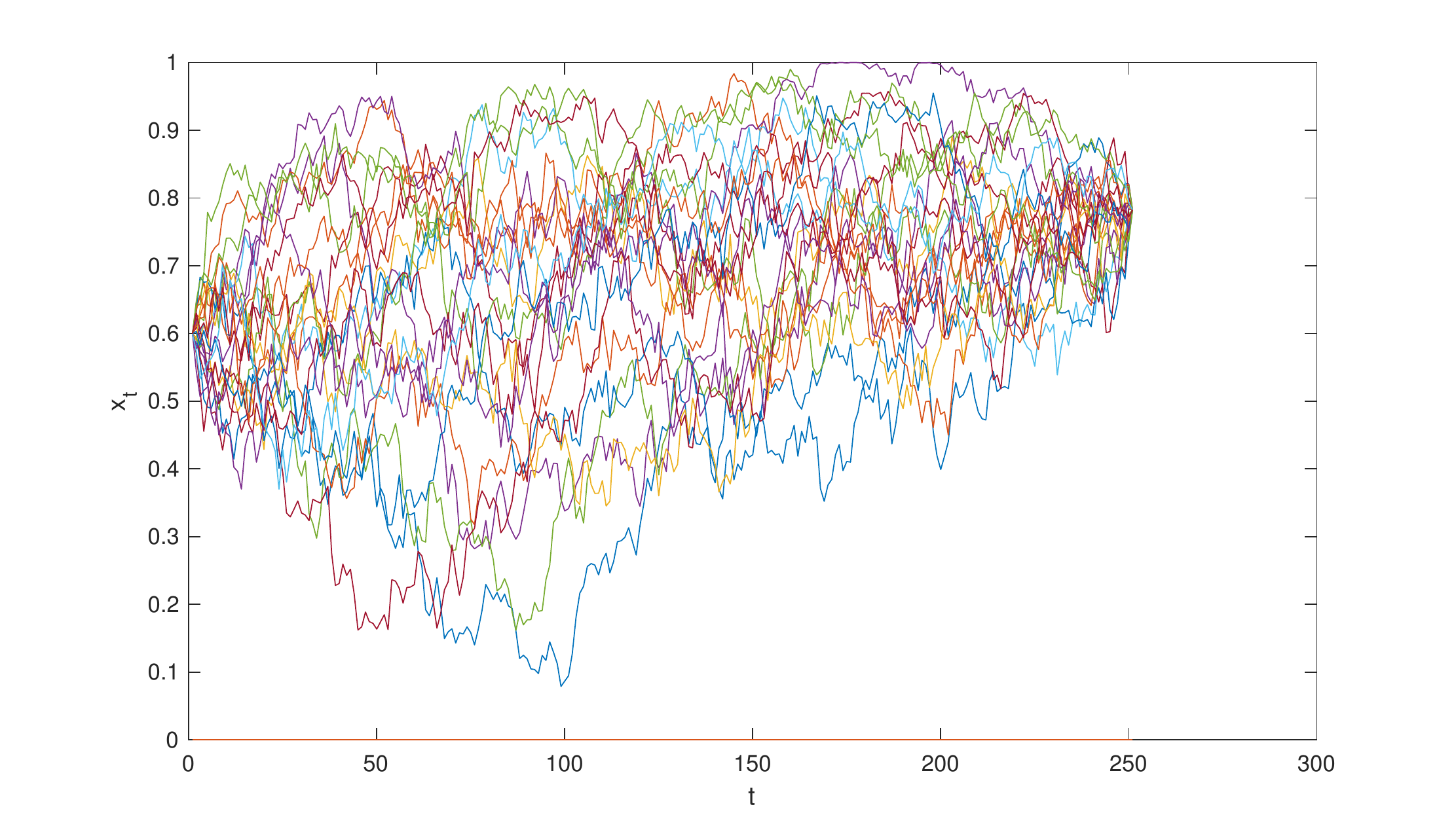}
 
    \caption{{\bf{Fig.4}} Diffusion bridge of accepted paths   }
  \label{fig:322}
\end{figure}

\subsubsection{The $2$-dimensional Wright-Fisher diffusion   with selection.}
In this section we investigate the performance our rejection sampling method in higher dimensions by considering the  $2-$dimensional Wright-Fisher diffusion with selection. Referring to  the work of Aurell et.al (2019)\cite{aurell2019multilocus} where the 2-dimensional Wright-Fisher diffusion with selection only is given by the stochastic  differential equation,
\begin{eqnarray}
       \label{eqn4350}
  \begin{cases}
      &dx^{(1)}_t=hx^{(1)}_t(1-x^{(1)}_t)x^{(2)}_tdt+\sqrt{x^{(1)}_t(1-x^{(1)}_t) } dW^{(1)}_t \\
            & dx^{(2)}_t=hx^{(2)}_t(1-x^{(2)}_t)x^{(1)}_tdt+\sqrt{x^{(2)}_t(1-x^{(2)}_t) } dW^{(2)}_t 
\end{cases} 
   ,
           \end{eqnarray} 
       where $h$ is the selection coefficient, $W^{(i)}_t\,\,\, i=1,2$ are $1$dimensional independent Wiener process  and  initial conditions $(x^{(1)}_0, x^{(2)}_0)=(x, y).$ 
       
       Similar to the $1-$dimensional case, we consider the transformation of each component of  the 2-dimensional standard Wiener process by the sigmoid function equation  (\ref{sigmoid}). By Ito's Lemma, the derived processes  is given by the following SDE,      
\begin{eqnarray}
       \label{eqn_42}
  \begin{cases}
      & dy^{(1)}_t=\frac{1}{2}y^{(1)}_t(1-y^{(1)}_t)(1-2y^{(1)}_t)dt+{y^{(1)}_t(1-y^{(1)}_t) } dW^{(1)}_t\\
      & dy^{(2)}_t=\frac{1}{2}y^{(2)}_t(1-y^{(2)}_t)(1-2y^{(2)}_t)dt+{y^{(2)}_t(1-y^{(2)}_t) } dW^{(2)}_t,
\end{cases} \end{eqnarray}  
with initial conditions  $(y^{(1)}_0, y^{(2)}_0)=(X, Y).$\\      
 Using the transformation theorem of random variables, it can be shown that the transition density of the process (\ref{eqn_42})  is given by   
 
 \begin{eqnarray}
\label{eqn_43 }
  P_1(t_0,y^{(1)}_0,y^{(2)}_0; t, y^{(1)}_t, y^{(2)}_t)=\begin{cases}
      
&\frac{1}{{2\pi (t-t_0)\sqrt{(1-\rho^2)}}}\times\frac{1}{y^{(1)}_ty^{(2)}_t(1-y^{(1)}_t)(1-y^{(2)}_t)} \times \\  
 & \exp \Big\{\frac{ (\beta_t^{(1)}-\beta_0^{(1)})^2+(\beta_t^{(2)}-\beta_0^{(2)})^2- 2(\beta_t^{(1)}-\beta_0^{(1)})(\beta_t^{(2)}-\beta_0^{(2)})\rho}{2(t-t_0)(1-\rho^2)}\Big\}
   \end{cases} \end{eqnarray} 
    where $\beta^{(i)}_t=\log{(\frac{y^{(i)}_t}{1-y^{(i)}_t})}\,\,\, i=1,2$  and $\rho$ is the correlation coefficient.  
    
    Now, considering  the process given by equation (\ref{eqn4350}),   under the multivariate transformation $ y^{(i)}=\frac{1}{1+e^{-\sin^{-1}(2x^{(i)}-1)}} \,\,\,i=1,2, $
    we obtain a SDE satisfied by the transformed Wright-Fisher process. 
    \begin{eqnarray}
    \label{eqn_433}
    \begin{cases}
      & dy^{(1)}_t=\frac{1}{2}y^{(1)}_t(1-y^{(1)}_t)(1-2y^{(1)}_t-\tan{\beta_t^{(1)}}+\\ &\frac{h}{2}(1+\sin {\beta_t^{(2)}}) \cos {\beta_t^{(1)}})dt+{y^{(1)}_t(1-y^{(1)}_t) } dW^{(1)}_t\\ 
      & dy^{(2)}_t=\frac{1}{2}y^{(2)}_t(1-y^{(2)}_t)(1-2y^{(2)}_t-\tan{\beta_t^{(2)}}+\\ &\frac{h}{2}(1+\sin {\beta_t^{(1)}}) \cos {\beta_t^{(2)}})dt+{y^{(2)}_t(1-y^{(2)}_t) } dW^{(2)}_t.
    \end{cases} 
    \end{eqnarray} 
     
  Denote by $P_1$ and $P_2$ the transition densities of the processes (\ref{eqn_433}) and (\ref{eqn_42})  respectively. Let $V(y_1,y_2,t)=\frac{P_2}{P_1}(y_1,y_2,t)$.  By theorem $3.1$, the ratio of transition densities $V$  satisfies the following partial differential equation,
    \begin{eqnarray}
\label{eqn_44}
\frac{\partial V}{\partial{t}}&=&a(y_1,y_2,t)V+b_1(y_1,y_2,t)\frac{\partial V}{\partial y_1}+b_2(y_1,y_2,t)\frac{\partial V}{\partial y_2}  +c_1(y_1,y_2)\frac{\partial^2 V}{\partial{y_1^2}}+c_2(y_1,y_2)\frac{\partial^2 V}{\partial{y_2^2}},\,\,\, t>0,\nonumber \\
\end{eqnarray} where
$$a(y_1,y_2,t)=\frac{1}{2}
 \begin{cases}\sec^2\beta_1+ \sec^2\beta_2+\frac{h}{2} (\sin \beta_1+\sin \beta_2+2\sin \beta_1\sin \beta_2)+\\
 \frac{\beta_1-\beta^{(1)}_0 -2(\beta_2-\beta^{(2)}_0)\rho}{t(1-\rho^2)}(\frac{h}{2}(1+\sin \beta_2)\cos \beta_1 -\tan \beta_1)+\\
  \frac{\beta_2-\beta^{(2)}_0 -2(\beta_1-\beta^{(1)}_0)\rho}{t(1-\rho^2)}(\frac{h}{2}(1+\sin \beta_1)\cos \beta_2 -\tan \beta_2), \end{cases} $$
 $$b_i(y_1,y_2,t)= \frac{1}{2}y_i(1-y_i)(1-2y_i+\frac{2[\beta^{(i)}_0-\beta_i +2(\beta_j-\beta^{(j)}_0)\rho]}{t(1-\rho^2)}+\tan \beta_i -\frac{h}{2}(1+\sin \beta_j)\cos \beta_i )  $$
 for $i\not=j, \,\,\,i,j=1,2.$
   
  $$c_i(y_1,y_2)=\frac{y_i^2(1-y_i)^2}{2}\,\,\, i=1,2.$$
  Due to the non existence of numerical methods for solving a general two dimensional parabolic partial differential equation Mohanty $\&$ Jain (1996)\cite{mohanty1996high},  the following change of coordinates is employed to obtain a PDE with non variable diffusion coefficients.
   
Let $x=x(y_1,y_2)=\beta_1$  and $y=y(y_1,y_2)=\beta_2$. Under this change of coordinates, equation (\ref{eqn_44}) is transformed into 
 
 \begin{eqnarray}
\label{eqn_45}
\frac{\partial V}{\partial{t}}&=&2\epsilon(x,y,t)V+c(x,y,t)\frac{\partial V}{\partial x}+d(x,y,t)\frac{\partial V}{\partial y}  +\frac{1}{2}\frac{\partial^2 V}{\partial{x^2}}+\frac{1}{2}\frac{\partial^2 V}{\partial{x^2}},\,\,\, t>0,\nonumber \\
\end{eqnarray}   
where
$$\epsilon(x,y,t)=\frac{1}{4}
 \begin{cases}\sec^2x+ \sec^2y+\frac{h}{2} (\sin x+\sin y+2\sin x\sin y)+\\
 \frac{x-\beta^{(1)}_0 -2(y-\beta^{(2)}_0)\rho}{t(1-\rho^2)}(\frac{h}{2}(1+\sin y)\cos x -\tan x)+\\
  \frac{y-\beta^{(2)}_0 -2(x-\beta^{(1)}_0)\rho}{t(1-\rho^2)}(\frac{h}{2}(1+\sin x)\cos y -\tan y), \end{cases} $$
 $$c(x,y,t)= \frac{\beta^{(1)}_0-x +2(y-\beta^{(2)}_0)\rho}{t(1-\rho^2)}+\frac{1}{2}\tan x-\frac{h}{4}(1+\sin y)\cos x, $$   
 $$d(x,y,t)= \frac{\beta^{(2)}_0-y +2(x-\beta^{(1)}_0)\rho}{t(1-\rho^2)}+\frac{1}{2}\tan y-\frac{h}{4}(1+\sin x)\cos y.  $$
 Following Karaa (2009)\cite{karaa2009high} for the numerical solution to (\ref{eqn_45}), define  difference operators in the coordinate directions $x$ and $y$ as follows,
 $$ A_x=-\big[ \alpha +\frac{\triangle x^2}{12}(\frac{c^2}{\alpha }-\epsilon -2\delta_x c)\big]\delta_x^2+\big[ c+\frac{\triangle x^2}{12}(\delta^2_x c-\frac{c}{\alpha }\delta_x c+2\delta_x \epsilon  -\frac{c\epsilon}{\alpha})\big]\delta_x+  \big[ \epsilon +\frac{\triangle x^2}{12}(\delta^2_x \epsilon-\frac{c}{\alpha }\delta_x \epsilon \big] ,$$
 $$L_x=\big[ 1+\frac{\triangle x^2}{12}(\delta^2_x -\frac{c}{\alpha }\delta_x) \big]$$ and 
 $$ A_y=-\big[ \alpha +\frac{\triangle y^2}{12}(\frac{d^2}{\alpha }-\epsilon -2\delta_y d)\big]\delta_y^2+\big[ d+\frac{\triangle y^2}{12}(\delta^2_y d-\frac{d}{\alpha }\delta_y+2\delta_y \epsilon  -\frac{d\epsilon}{\alpha})\big]\delta_y+  \big[ \epsilon +\frac{\triangle y^2}{12}(\delta^2_y \epsilon-\frac{d}{\alpha }\delta_y \epsilon \big] ,$$ 
 $$L_y=\big[ 1+\frac{\triangle y^2}{12}(\delta^2_y -\frac{d}{\alpha }\delta_y \big].$$   
 Here, $\alpha$ is the diffusion coefficient  in the parabolic PDE  which in this case  is equal to $\frac{1}{2}$.  $\delta^2_{-}$  and $\delta_{-}$ are the central difference operators  for the approximation of the  second and first derivatives in the respective directions, $\triangle x$ and $\triangle y$ are respectively the spatial step length in the $x$ and $y$ direction. \\In addition, denote the  length of time step by $\triangle t$.

 Under the assumption $\triangle t \leq min(\triangle x, \triangle y)$,  Karaa (2009)\cite{karaa2009high} showed that, the Cranck-Nicolson scheme for the parabolic PDE  (\ref{eqn_45}) is given by
 \begin{eqnarray}
\label{eqn_46}
(L_x+\frac{\triangle t}{2}A_x)(L_y+\frac{\triangle t}{2}A_y)V^{n+1}_{i,j}=(L_x-\frac{\triangle t}{2}A_x)(L_y-\frac{\triangle t}{2}A_y)V^{n}_{i,j},
\end {eqnarray}
 where $V_{i,j}$ is the approximate solution to $V$ at the point $(x_i,y_j)$ at time step $n$. 
 In implementing the numerical  scheme  (\ref{eqn_46}), it is split into  the alternating direction implicit (ADI)  scheme  
  \begin{eqnarray}
\label{eqn_47}
(L_x+\frac{\triangle t}{2}A_x)V^{\ast}=(L_x-\frac{\triangle t}{2}A_x)(L_y-\frac{\triangle t}{2}A_y)V^{n}_{i,j},\\
(L_y+\frac{\triangle t}{2}A_y)V^{n+1}=V_{i,j}^{\ast}
\end {eqnarray} where $V^{\ast}$ is the approximate solution at an intermediate time step between $t=n$ and $t=n+1.$
Utilising the the above ADI scheme, the normalised numerical solution $V$  to the PDE at different time steps are shown in Figures $5$.
 
  \begin{figure}[h]
  \centering
   \includegraphics[ width=.43\textwidth]{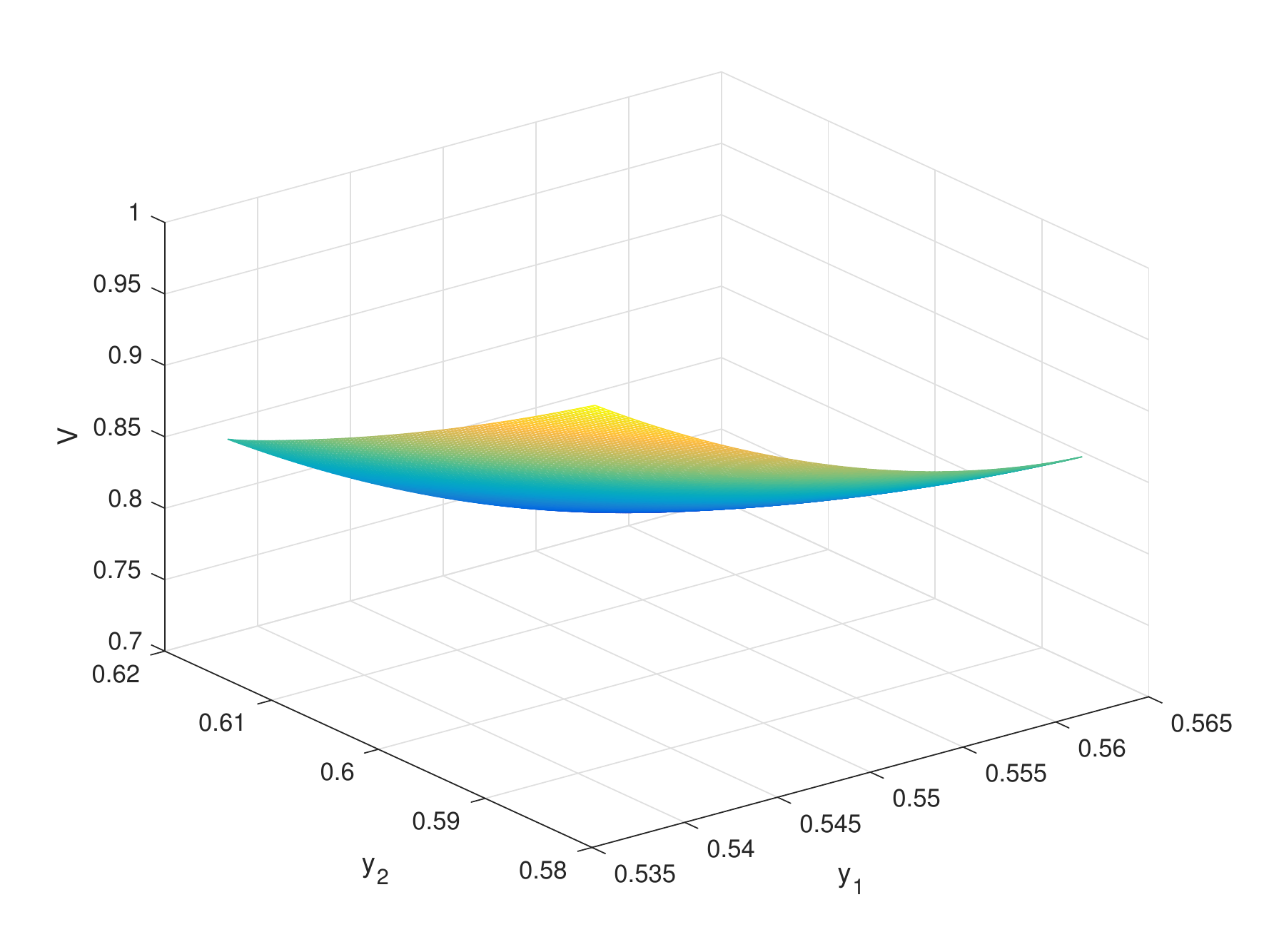} 
    \includegraphics[ width=.43\textwidth]{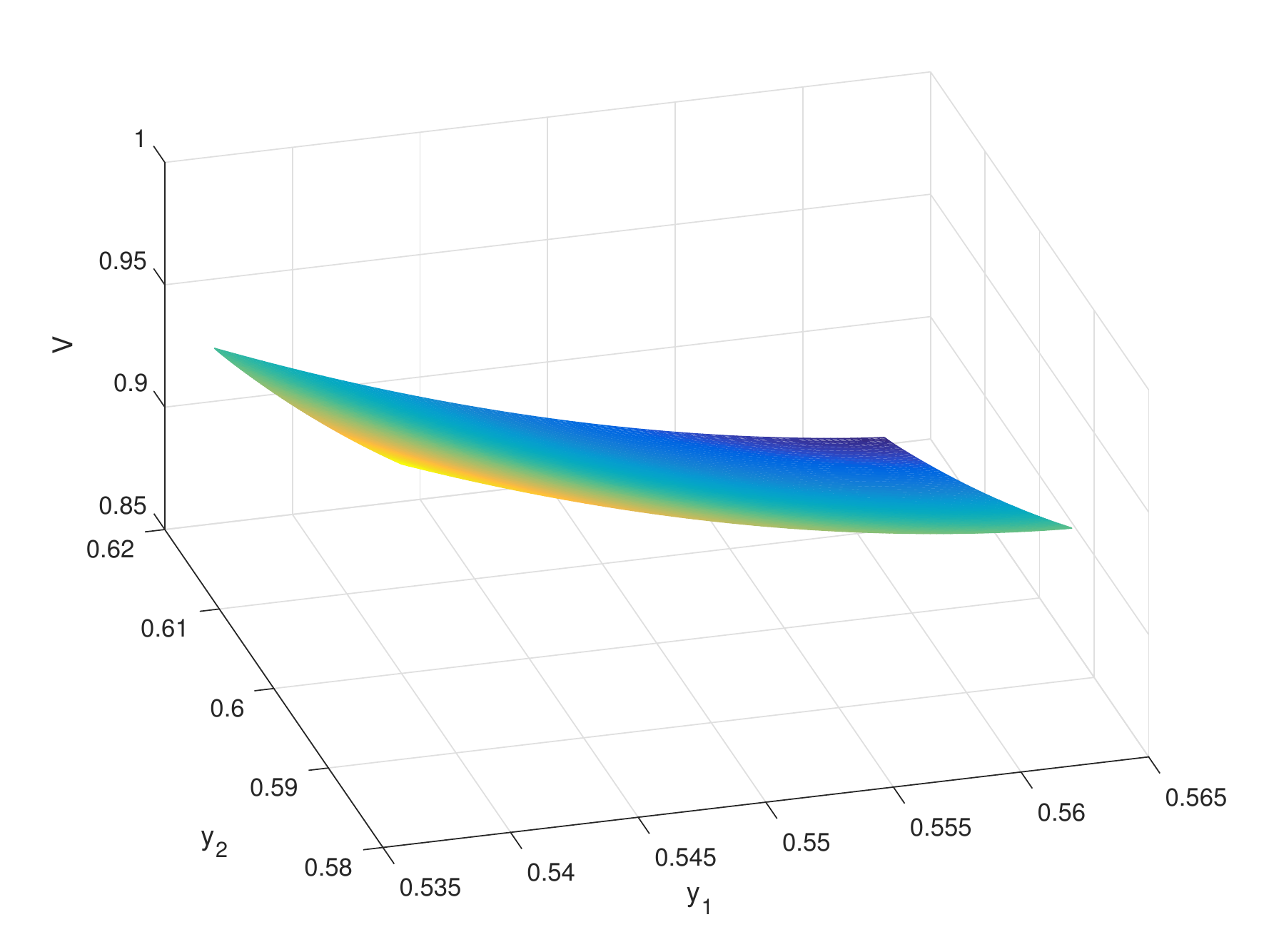}
  \includegraphics[ width=.43\textwidth]{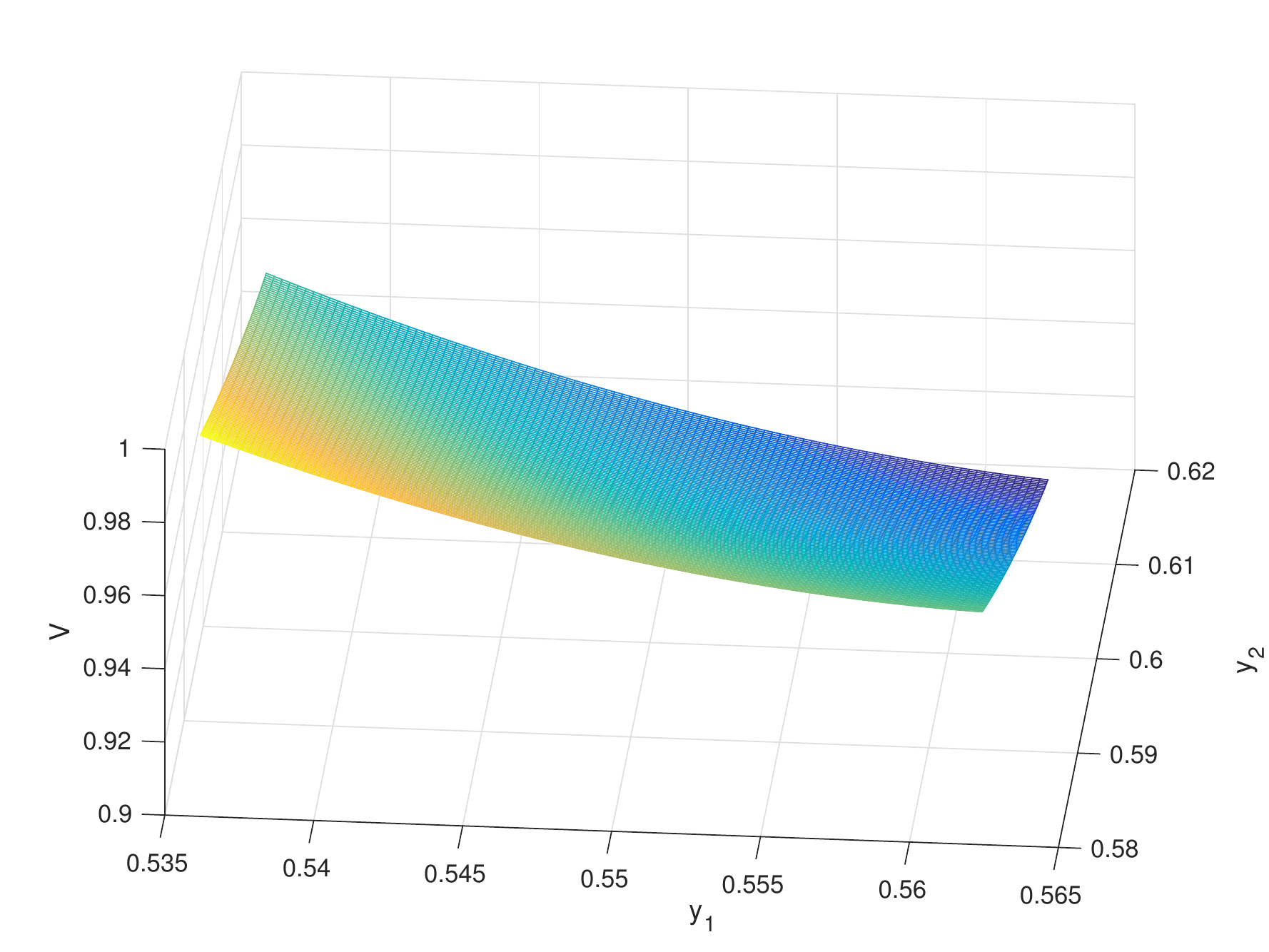}
    \includegraphics[ width=.43\textwidth]{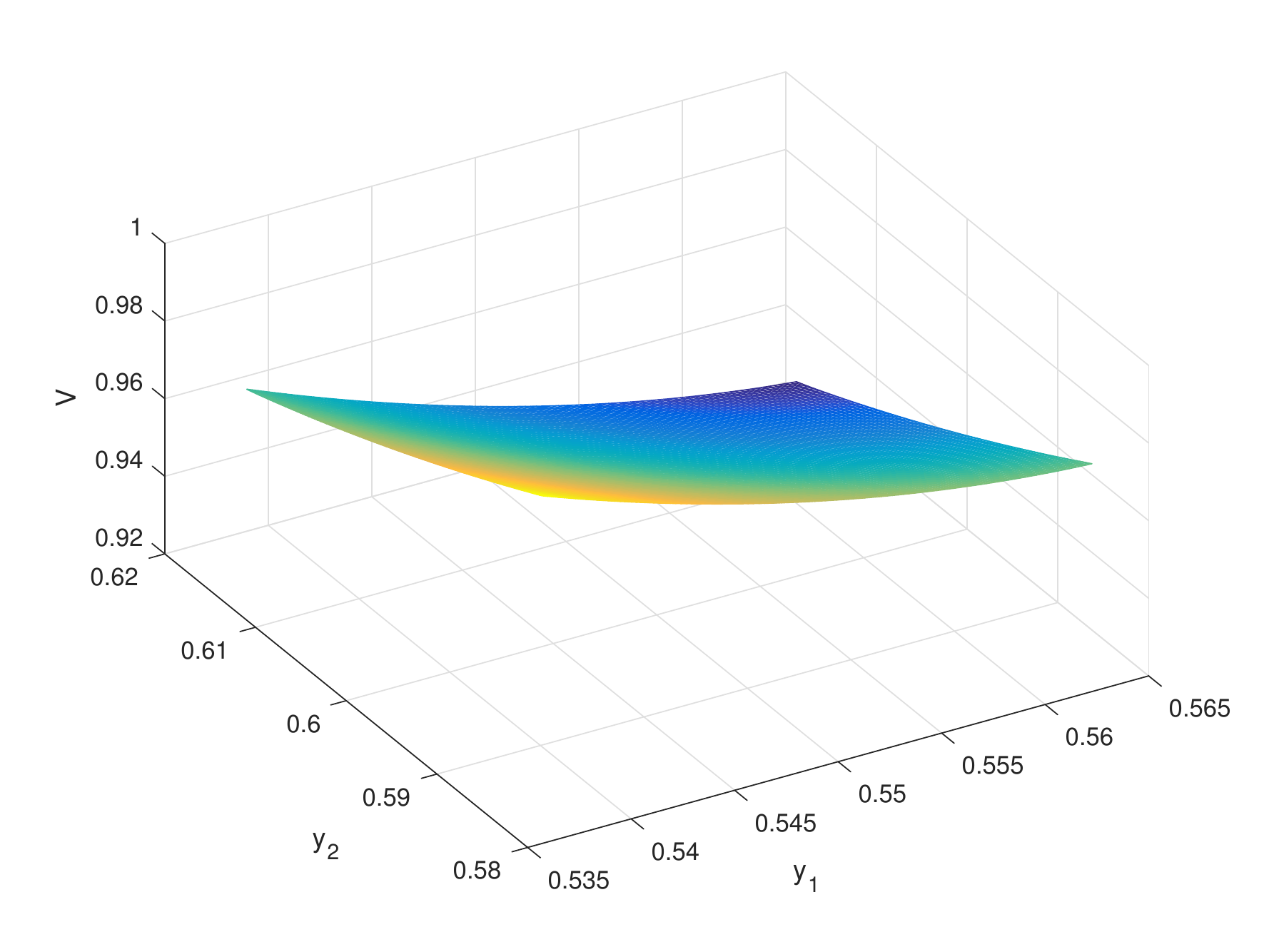}
     \caption{{\bf{Fig.5}} Normalised numerical solution $V$ at time step $5,10,15,20$ respectively for $h=1$}
       \label{fig:350}
\end{figure}

The numerical solution displays a decreasing behaviour for increasing values on both axes which is rather expected since for a fixed value of $h$, there is an increasing divergence between the proposal process and the target process.  A similar  property of the solution is further displayed by the solution for  increasing values of the parameter $h$ shown in the appendix, which is explained by the same reason given in the previous line.  Thus, a higher acceptance probability is expected for a process which spends a significant  time near the left lower corner of the square $(0,1)\times(0,1)$ and counter behaviour  at the right top corner, see  second plot in Figure $6$. This is explained by  fact that the proposal process best approximates  the target process for values near the zero point.

 \begin{figure}[h]
 
    \centering
  \includegraphics[ width=0.43\textwidth]{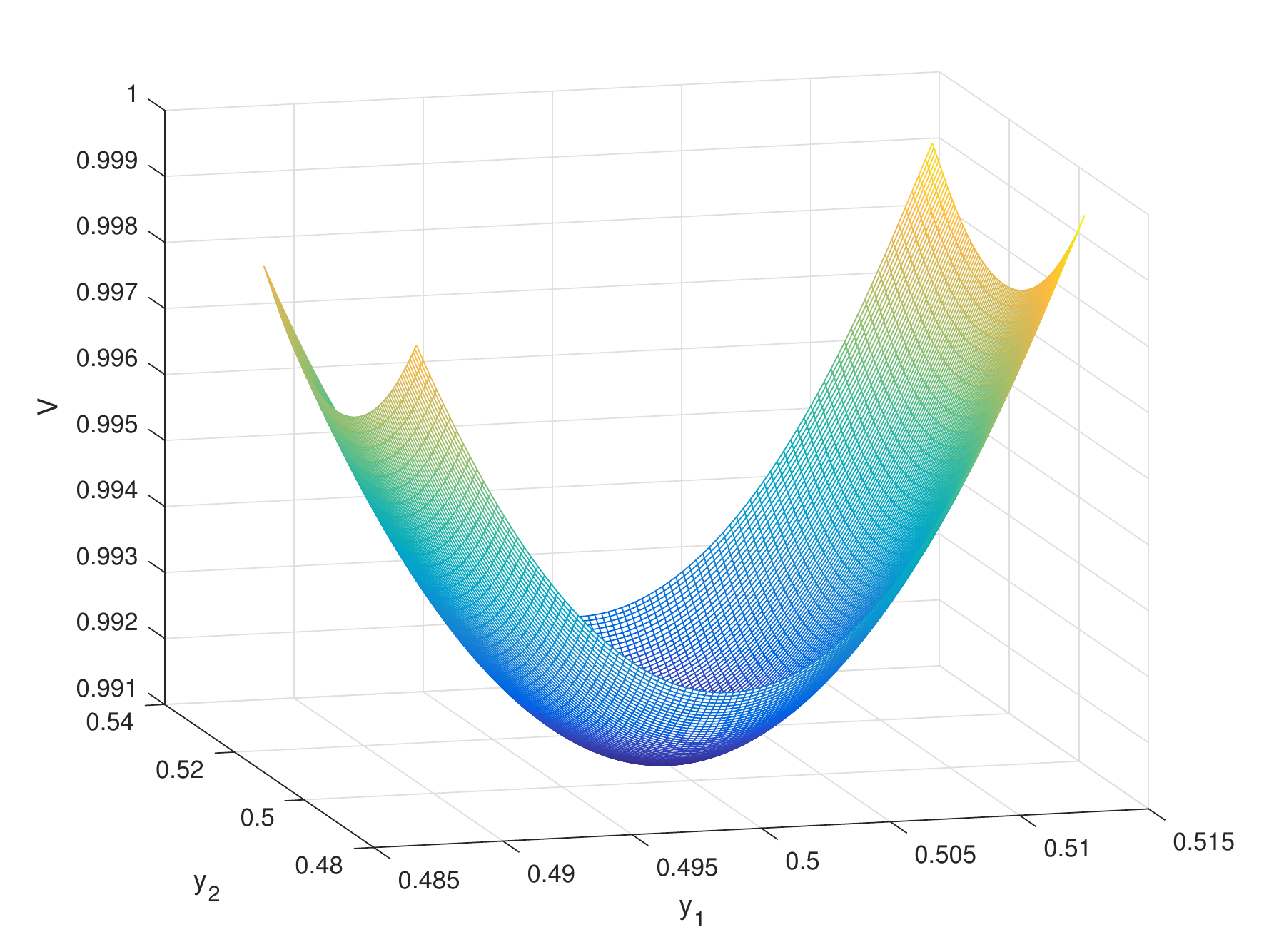}
   \includegraphics[ width=.43\textwidth]{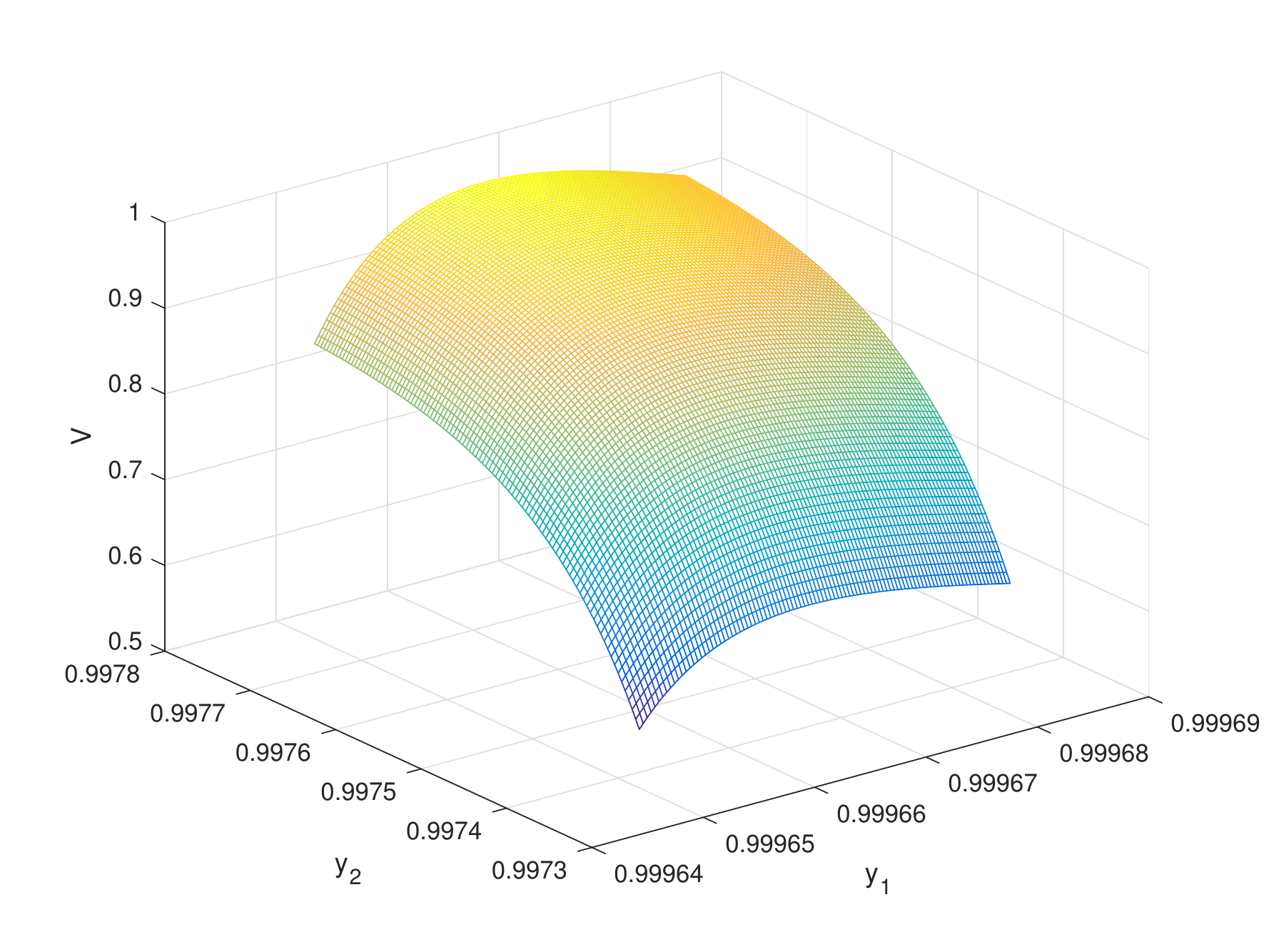}
        \caption{{\bf{Fig.6}} Normalised numerical  solution $V$ at time step  $20$  for $h=0$}
  \label{fig:321}
\end{figure}
However, we observe a slightly different behaviour in the numerical solution for $h=0$, see the \\left hand plot in Figure $6$. Here,  the proposal process differs from the target process by only terms $\tan(\beta^{(i)})$ which influences the properties of the solution. 
\vspace{5cm}

  \subsection{Unbounded Drift: The  Ornstein -Uhlenbeck  Process  } 
  Here,  we implement the our method   on particular processes  with  known explicit forms of  transition densities and  take a comparison between the approximated ratio values using the PDE and the theoretical values for processes.  Particularly, we consider the Ornstein -Uhlenbeck  Process  whose drift is unbounded and hence its  transition density Kusuoka (2017)\cite{kusuoka2017continuity}. We consider the Wiener process  as the reference process in our rejection sampling.  \\ 
    Consider an  Ornstein -Uhlenbeck  process  given the stochastic differential equation,
   \begin{eqnarray}
\label{eqn8}
dX_{t}&=&-\beta X_tdt+\sigma dW_t, \,\,\, X_0=x_0, 0\leq t \leq T,\,\,\, \beta>0. 
\end {eqnarray} 
 Suppose we require to sample from the   conditional densities of the O-U process at different times  using   Wiener process as the reference process.
 We thus consider the Wiener process given by the SDE
   \begin{eqnarray}
\label{eqn8a}
dX_{t}&=&\sigma dW_t, \,\,\, X_0=x_0, 0\leq t \leq T. 
\end {eqnarray}  
 The Fokker-Planck equations associated with the equations (\ref{eqn8a}) and (\ref{eqn8}) are respectively  given  by
\begin{eqnarray}
\label{eqn9}
\frac{\partial}{\partial{t}}{P_1(x, t/x_0,t_0)}&=&\frac{1}{2}\frac{\partial^2}{\partial{x^2}}{[\sigma^2 P_1(x, t/x_0,t_0)]},\\
\label{eqn9b}
\frac{\partial}{\partial{t}}{P_2(x, t/x_0,t_0)}&=&\frac{\partial}{\partial{x}}[\beta xP_2(x, t/x_0,t_0)]+\frac{1}{2}\frac{\partial^2}{\partial{x^2}}{[\sigma^2 P_2(x, t/x_0,t_0)]},
\end{eqnarray} 
 with initial condition $ P_i(x, t_0/x_0,t_0)=\delta(x-x_0)\,\,\,\,i=1,2$.\\ 
 However, by Ito calculus,  $P_2(x, t/x_0,t_0)$  is known to be  a Gaussian density with mean\\ $\mu_2(t)=x_0e ^{-\beta t}$  and variance $\sigma_2^2(t)=\frac{\sigma^2}{2\beta}(1-e^{-2\beta t})$. \\Similarly, $P_1(x, t/x_0,t_0)\in N( \mu_1(t),\sigma_1^2(t))$, where  $\mu_1(t)=x_0$   and $\sigma_1^2(t)=\sigma^2 t$.
 {\propn
The quotient of transition densities  $\frac{P_2}{P_1}(x, t/x_0,t_0)$ satisfies the parabolic partial differential equation (\ref{eqn6})}.

\begin{proof}
  $\sigma(x)=\sigma^2$, $S_1(x)=0$, $S_2(x)=-\beta x$, $S_1'(x)=0$ and $S'_2(x)=-\beta$, 
$f(x,t)=\frac{\partial }{\partial x}\log{P_1}$  implying   $f(x,t)= -\frac{x-\mu_1}{\sigma_1^2}$.
\begin{eqnarray}
a(x,t)& = &  S'_1(x)-S'_2(x)+[S_1(x)-S_2(x)]f(x,t), \nonumber 
\end{eqnarray}
implying that \begin{eqnarray}
 a(x,t)& = & \frac{\beta}{\sigma_1^2}\Big\{\sigma_1^2 -{x(x-\mu_1)}\Big\}, 
\end{eqnarray}\begin{eqnarray}
b(x,t) & = & \sigma(x)f(x,t)+\sigma'(x)-S_2(x), \nonumber   \\
 & = & (\beta-\frac{\sigma^2}{\sigma_1^2} )x+\frac{\sigma^2\mu_1}{\sigma_1^2}
\end{eqnarray}
and  
\begin{equation}
c(x,t)=\frac{\sigma^2}{2}.
\end{equation}

\begin{equation}
\label{ }
\frac{\partial}{\partial{x}}{\big[\frac{P_2}{P_1}\big]}=\frac{P_2}{P_1}\Big\{ (\frac{1}{\sigma_1^2} - \frac{1}{\sigma_2^2})x+\frac{\mu_2}{\sigma_2^2}- \frac{\mu_1}{\sigma_1^2}\Big\}
\end{equation}
Let $A=(\frac{1}{\sigma_1^2} - \frac{1}{\sigma_2^2})x+\frac{\mu_2}{\sigma_2^2}- \frac{\mu_1}{\sigma_1^2}.$ Then, 

\begin{equation}
\frac{\partial^2}{\partial{x^2}}{\big[\frac{P_2}{P_1}\big]}=\frac{P_2}{P_1}\Big\{ A^2+\frac{1}{\sigma_1^2} - \frac{1}{\sigma_2^2}\Big\}.
\end{equation}
Considering the left hand term  of  equation \ref{eqn6},
\begin{equation}
\label{A2}
\frac{\partial}{\partial{t}}{\big[\frac{P_2}{P_1}\big]}=\frac{P_2}{P_1}\Big\{ \frac{\sigma_1'}{\sigma_1} - \frac{\sigma_2'}{\sigma_2} +(\frac{1}{\sigma_1}x-\frac{\mu_1}{\sigma_1})(\frac{-\sigma_1'}{\sigma_1^2}x+\frac{\mu'_1}{\sigma_1}+\frac{-\sigma_1'}{\sigma_1^2}\mu_1) - (\frac{1}{\sigma_2}x-\frac{\mu_2}{\sigma_2})(\frac{-\sigma_2'}{\sigma_2^2}x+\frac{\mu'_2}{\sigma_2}+\frac{-\sigma_2'}{\sigma_2^2}\mu_2   )   \Big\}.
\end{equation}
Let 
\begin{equation}
\label{A3}
h(x,t)=a(x,t)+b(x,t).A+c(x,t)(A^2+\frac{1}{\sigma_1^2} - \frac{1}{\sigma_1^2}).
\end{equation}
To check whether the partial differential equation ( \ref{eqn6}) is satisfied,  it is enough to check whether 
term in the right hand  parenthesis of  equation (\ref{A2}) is equal to equation (\ref{A3}) for all $(x,t)$.\\
However,   equations (\ref{A2}) and   (\ref{A3}) are quadratic polynomials in $x$ with time dependent  coefficients. \\From the laws of algebra, any two polynomials  are equal if and only if the coefficients are equal for all values t.
Thus, we check where the coefficients are equal and with some  manipulations of the involved terms equality is established. For instance, 
 comparing coefficients of the  quadratic terms, 
 from \ref{A2}, this is given by  
\begin{equation}
\label{t-derivative}
\frac{-\sigma_1'}{\sigma_1^3}+\frac{\sigma_2'}{\sigma_2^3}=-\frac{1}{2\sigma_1^4}\sigma^2+\frac{\sigma^2}{2\sigma_2^4}-\frac{\beta} {\sigma_2^2}.\end{equation}
The derivatives in (\ref{t-derivative}) are with respect to the time  variable $t$.

 Similarly, from (\ref{A3})  the coefficient of the quadratic term is given by 
\begin{eqnarray}
\frac{-\beta}{\sigma_1^2}+(\beta-\frac{\sigma^2}{\sigma_1^2})(\frac{1}{\sigma_1^2} - \frac{1}{\sigma_2^2})+\frac{\sigma^2}{2}(\frac{1}{\sigma_1^2} - \frac{1}{\sigma_2^2})^2 & = &
-\frac{1}{2\sigma_1^4}\sigma ^2+\frac{1}{2\sigma_2^4}\sigma^2-\frac{\beta}{\sigma_2^2}
\end{eqnarray}
Similarly the coefficients for the other  terms are also found to be  equal and thus the PDE( \ref{eqn6} ) is satisfied. \end{proof}

 Consequently, $\frac{P_2}{P_1}(x, t/x_0,t_0)$ evolves according to a second order linear  parabolic partial differential equation  whose coefficients  are given by 
  
$ a(x,t) = \frac{\beta}{\sigma_1^2}\Big\{\sigma_1^2 -{x(x-\mu_1)}\Big\}$, 
 $b(x,t)= (\beta-\frac{\sigma^2}{\sigma_1^2} )x+\frac{\sigma^2\mu_1}{\sigma_1^2}$
and $c(x,t)=\frac{\sigma^2}{2}$.

We further check for the bound on quotient $\frac{P_2}{P_1}$
since our  proposed rejection method  requires quotient $\frac{P_2}{P_1}$  to be bounded  for  $t\in [0,T].$ \\
 Consider the ratio of conditional densities $P_1$ and $P_2$,
 
\begin{eqnarray}
\label{ratio}
\frac{P_2}{P_1}(x,t/x_0,t_0) =\frac{\sigma_1}{\sigma_2}\exp{-\frac{1}{2}\Big\{{\frac{(x-x_0e^{-\beta t})^2}{\frac{\sigma^2}{2\beta}{(1-e^{-2\beta t})}}-{\frac{(x-x_0)^2}{\sigma^2 t}\Big\}}}}.
\end{eqnarray}
The  ratio (\ref{ratio}) is bounded above by a constant $C$ that depends on $\beta$. \\
That is, 
\begin{eqnarray}
\frac{\sigma_1}{\sigma_2}\exp{-\frac{1}{2}\Big\{{\frac{(x-x_0e^{-\beta t})^2}{\frac{\sigma^2}{2\beta}{(1-e^{-2\beta t})}}-{\frac{(x-x_0)^2}{\sigma^2 t}\Big\}}}}
\leq C(\beta).\end{eqnarray}

Hence for a fixed value of $\beta$ an upper  bound on ratio of conditional densities can be established.\\
Moreover,  such a bound for a ratio of  conditional transition densities with unknown   closed forms for diffusion  processes described by SDEs   can equally be established, see for 
instance, Downes (2008 )\cite{downes2009bounds} for a  detailed discussion of  bounds on transition densities of time-homogeneous diffusion  processes described by stochastic differential equations.\\
Following the work of Downes(2008)\cite{downes2009bounds} on bounds on transition densities,  it can be shown that  the upper bound $C(\beta)$ on $\frac{P_2}{P_1}$ in the case of the O-U process with Brownian motion as the reference process is given by

$$C(\beta)=e^{ \frac{\beta}{2}(x^2_{max}-x_0^2+T)}.$$
Now suppose the closed form of $ P_2$ is unknown but available is the SDE (\ref{eqn8}) for the stochastic process or equally the Fokker-Planck equation for $P_2$ equation  (\ref{eqn9b}) since to every Fokker-Planck equation we can associate a stochastic differential equation Gardiner (1985), chapter $5$. In addition, let the Wiener process (\ref{eqn8a}) be the proposal process in our rejection method.
Thus, we  require the solution to  parabolic 
 \begin{eqnarray}
\label{eqn12}
\frac{\partial}{\partial{t}}{V(x,t)}&=&\frac{\beta}{\sigma_1^2}\Big\{\sigma_1^2 -{x(x-\mu_1)}\Big\}V(x,t)+\nonumber \\ &&\Big\{(\beta-\frac{\sigma^2}{\sigma_1^2} )x+\frac{\sigma^2\mu_1}{\sigma_1^2}\Big \}\frac{\partial}{\partial{x}}{V(x,t)} +\frac{\sigma^2}{2}\frac{\partial^2}{\partial{x^2}}{V(x,t)},
\end{eqnarray}  
where $\mu_1=x_0$, $\sigma_1^2=\sigma^2t,$ subject to  initial  condition  $V(x,t_0)=1$. 
Noting that equation (\ref{eqn12}) is a second order parabolic equation, for solution we further require that we specify the boundary conditions. However, the boundary conditions depend on the nature of the boundaries  shared by both processes on the interval over which the processes are constrained over the time interval $ [0,T] $.  According to  Gardiner (1985)\cite{gardiner1985handbook}, the boundaries can take on any of the following forms;  reflecting, absorbing, periodic or prescribed boundaries.

In  this particular case,  the boundaries are finite considering the time period $[0,T]$ and  these are assumed  to be $x=x_{min}$ and $x=x_{max}$ and the boundary conditions are given by $V(x_{min},t)=1$ and $V(x_{max},t)=1$. \\
Considering the Cranck-Nicolson finite difference scheme for equation (\ref{eqn12})  as described in section 3.1,  he numerical solutions  for the PDE adjusted by $C(\beta)$ for 
 $\beta =25, 0.05, 0.005, 5\times10^{-20}$  with the corresponding simulated O-U processes are shown in the figure $1$

   \begin{figure}[h]
  \centering
  \includegraphics[ width=1\textwidth]{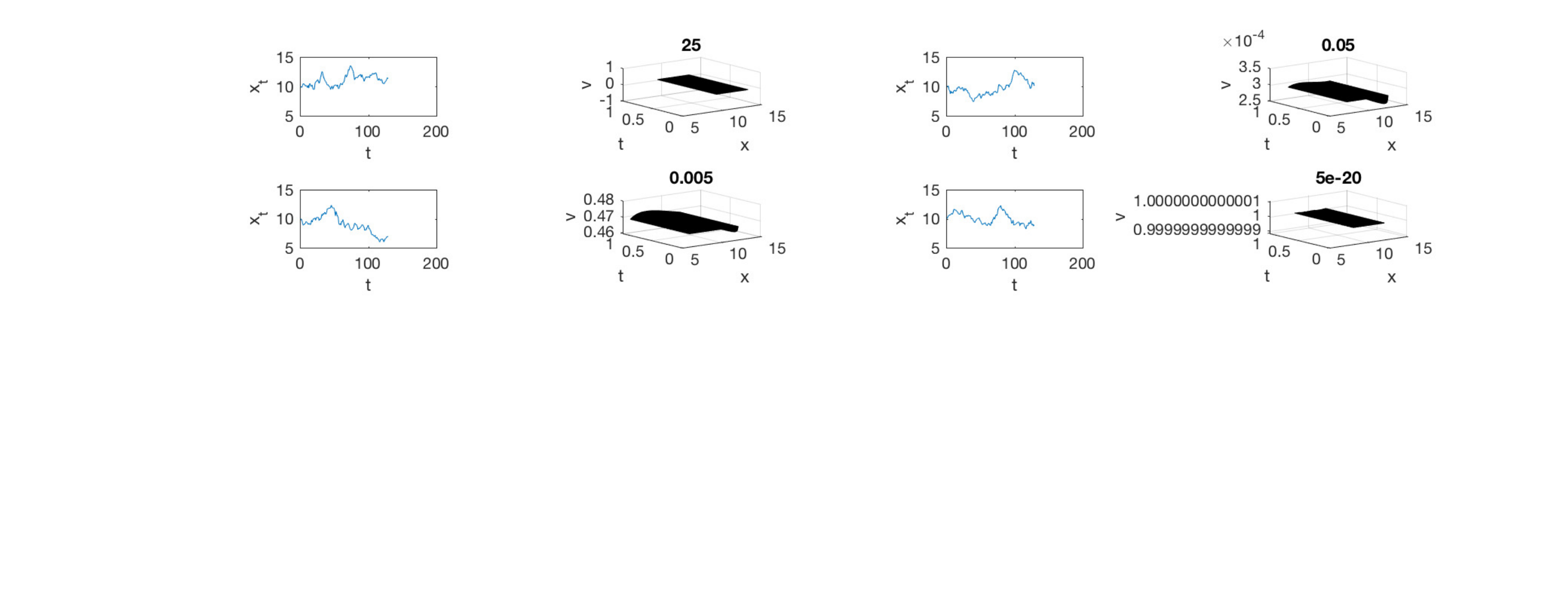}
  \vspace{-4cm}
 \caption{{\bf{Fig.7}} Normalised numerical solutions of the PDE for different $\beta$ values  with corresponding O-U processes}
  \label{fig:35}
\end{figure}
  \begin{figure}[h]
  \centering
  \includegraphics[ width=1\textwidth]{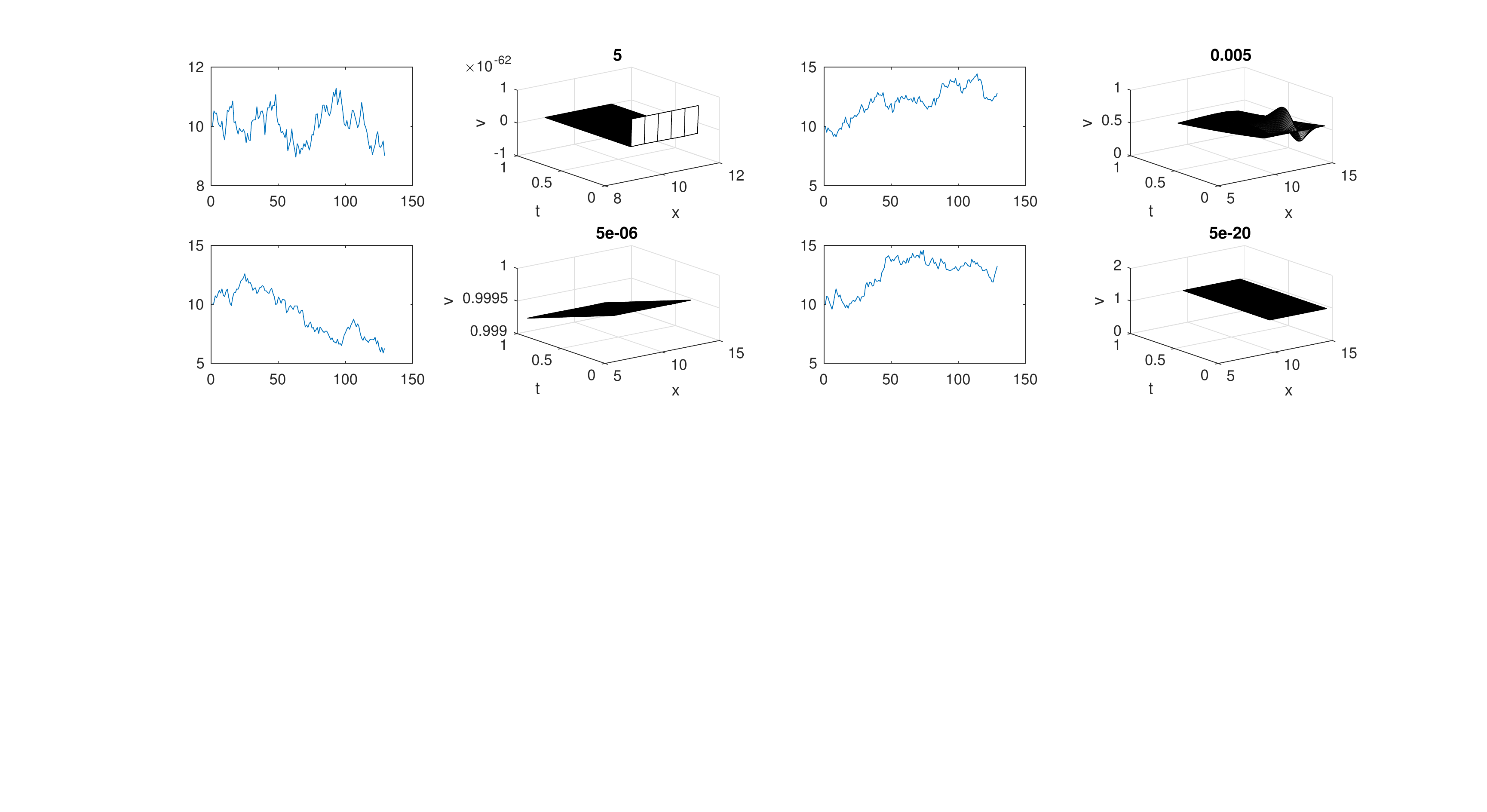}
   \vspace{-5cm}
    \caption{ {\bf{Fig.8}} Proposal discretized Brownian motion paths with corresponding normalised numerical solutions of the PDE for different $\beta$ values }
  \label{fig:36}
\end{figure}
 
\begin{figure}[h]
  \centering
  \includegraphics[ width=1\textwidth]{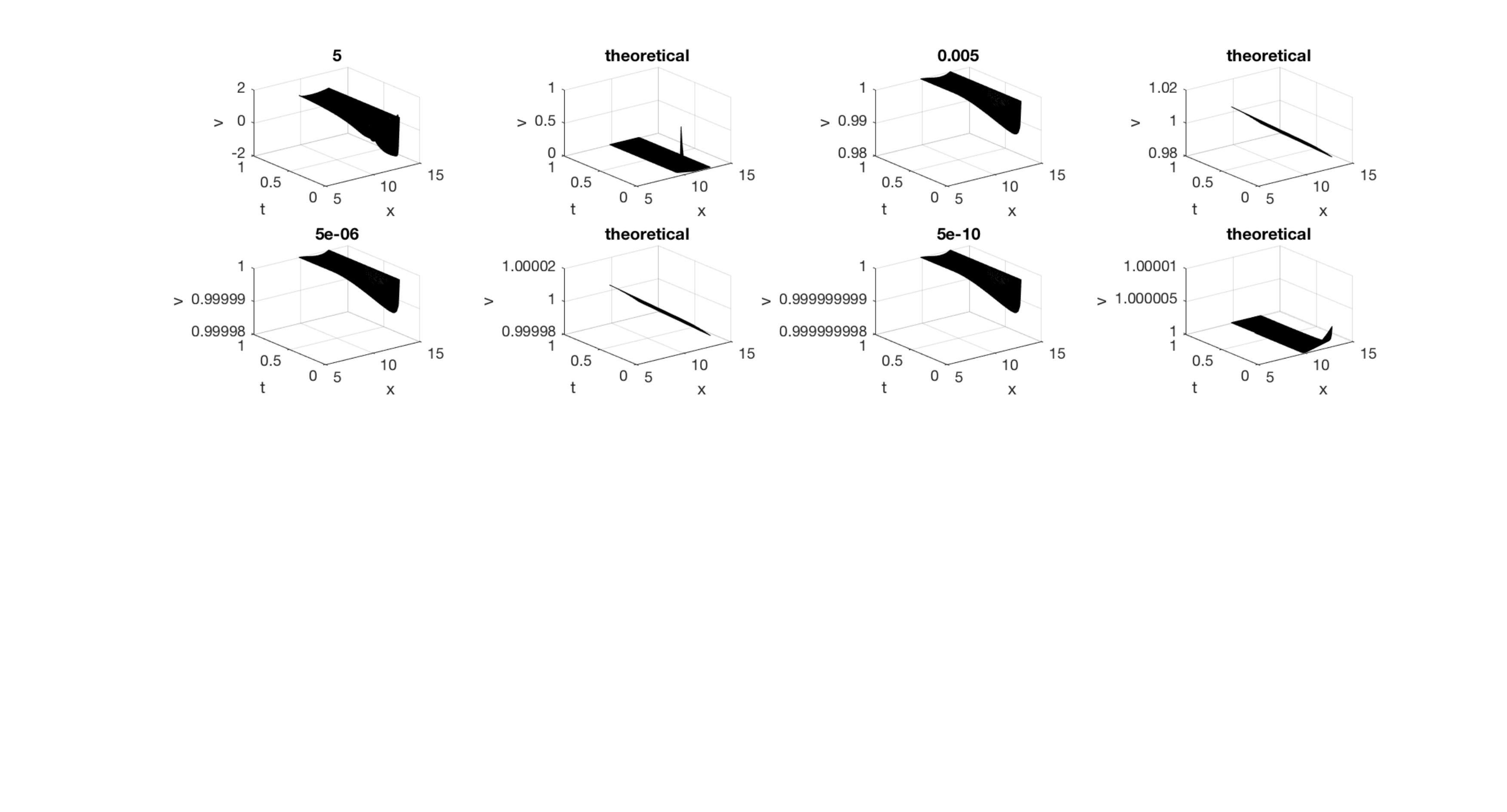}
  \vspace{-5cm}
 \caption{{\bf{Fig.9}} Normalised numerical solutions of the PDE for different $\beta$ values  with corresponding theoretical values}
  \label{fig:37}
\end{figure}

 The numerical solution to the PDE possesses  the expected behaviour of  $V$ tending to unity as $\beta$ tends to  zero. Observe that, 
 for large  values of $\beta$, $V$  tends to zero indicating that  $P_1$ is no longer enveloped in $P_2$ over the considered time interval.  This is  due to  the fact that  higher values of $\beta $ introduce a trend in the paths of the diffusion model (\ref{eqn8}) and thus candidate  paths with no trends can hardly be solutions.  This is the case  for $\beta=25$ in which the solutions to the SDE constitutes a trend. \\ Whereas $V$ tends to unity for very small $\beta$ values, for extremely large values of $\beta$ the  approximate solution to the PDE remains bounded. \\
 The proposed paths figure {\ref{fig:36}} are rejected in the case of $\beta=5$ for the reason  similar to the case of  $\beta=25$ mentioned in the previous paragraph.   However, for $\beta=0.005$  some  points on the path are rejected whereas others are accepted.  Thus,  the solution path can be constructed by interpolating  between the accepted points by means of a diffusion  bridge, a brownian bridge in this case.  We further note that for  the last two cases in  figure {\ref{fig:36}} , the proposed paths are accepted entirely and the solution at any other time  points  can similarly be obtained by Brownian bridges.\\
 Moreover,  the numerical solution to the PDE converges to the theoretical values as $\beta$  becomes small. That is, the absolute error in the numerical solution reduces is small for smaller values of $\beta$, see Figure  \ref{fig:37}. 

 \section{ Conclusion}
 In a nut shell therefore, we have proposed a new method of carrying out rejection sampling  using transition densities including those that are unknown in closed form. The  method has been tested on cases with known transition densities  and  the results  are consistent with theoretical values. Thus the method provides a possibility of simulating paths for  diffusions at discrete time points  using  associated  transition densities even when they are unknown in closed form.   We have  provided a multidimensional extension of the method in addition,  and we believe avails  an opportunity for use of   likelihood free inferential methods for numerous diffusion processes  such as  Approximately Bayesian  computational (ABC) methods. Since by our method sampling from the transition density of multidimensional Wright-Fisher diffusion at discrete time points is possible,   in the future, we embark on applying likelihood free methods in  determining the interaction graphical networks in allele frequency data using coupled  Wright-Fisher diffusion as the underlying model. 

 \newpage
\bibliographystyle{plain}
\bibliography{reference}
\section*{Appendix}
\renewcommand{\theequation}{A\thesection.\arabic{equation}}
 \begin{figure}[h]
  \centering
  \includegraphics[ width=.43\textwidth]{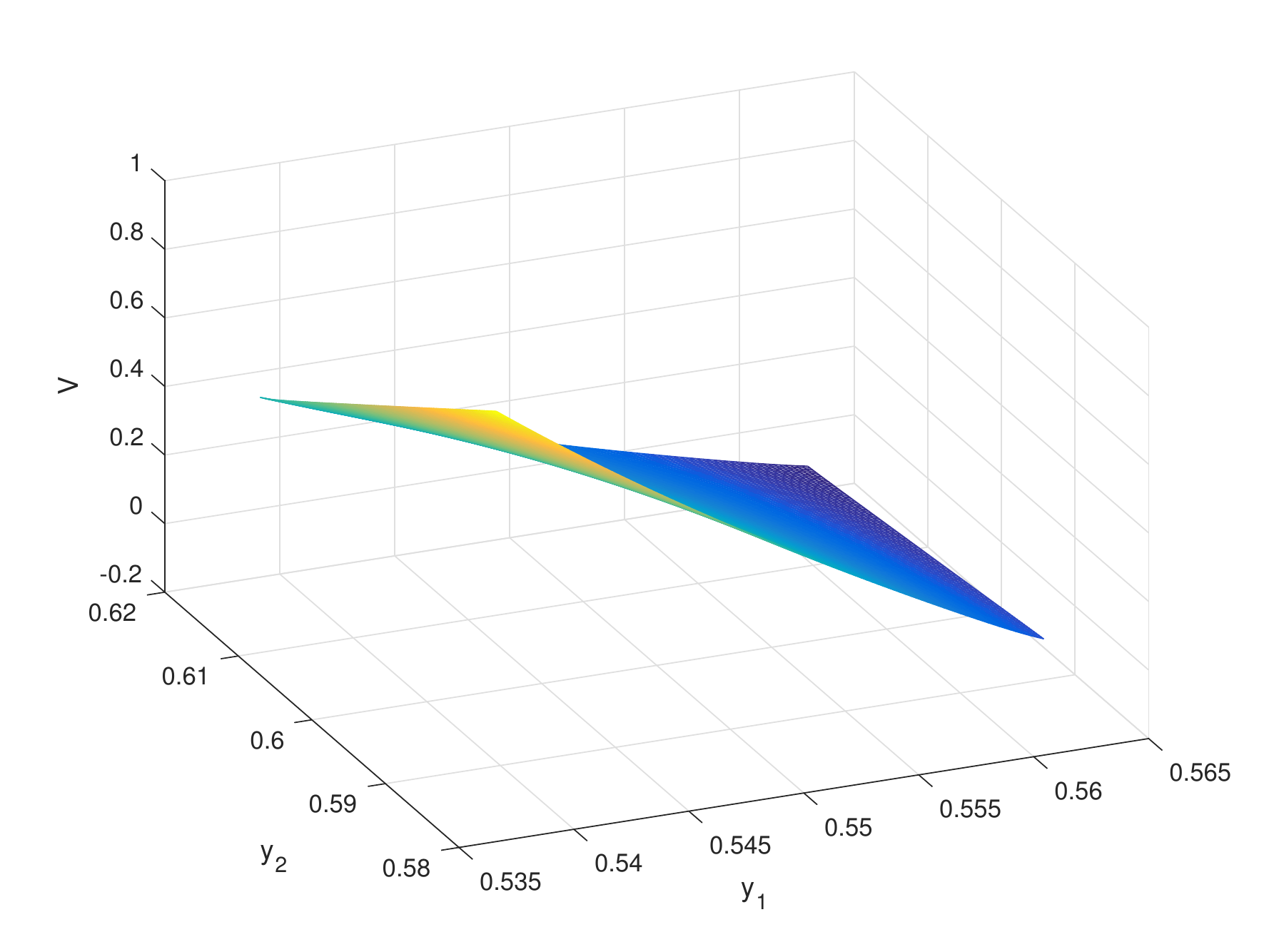}
  \includegraphics[ width=.43\textwidth]{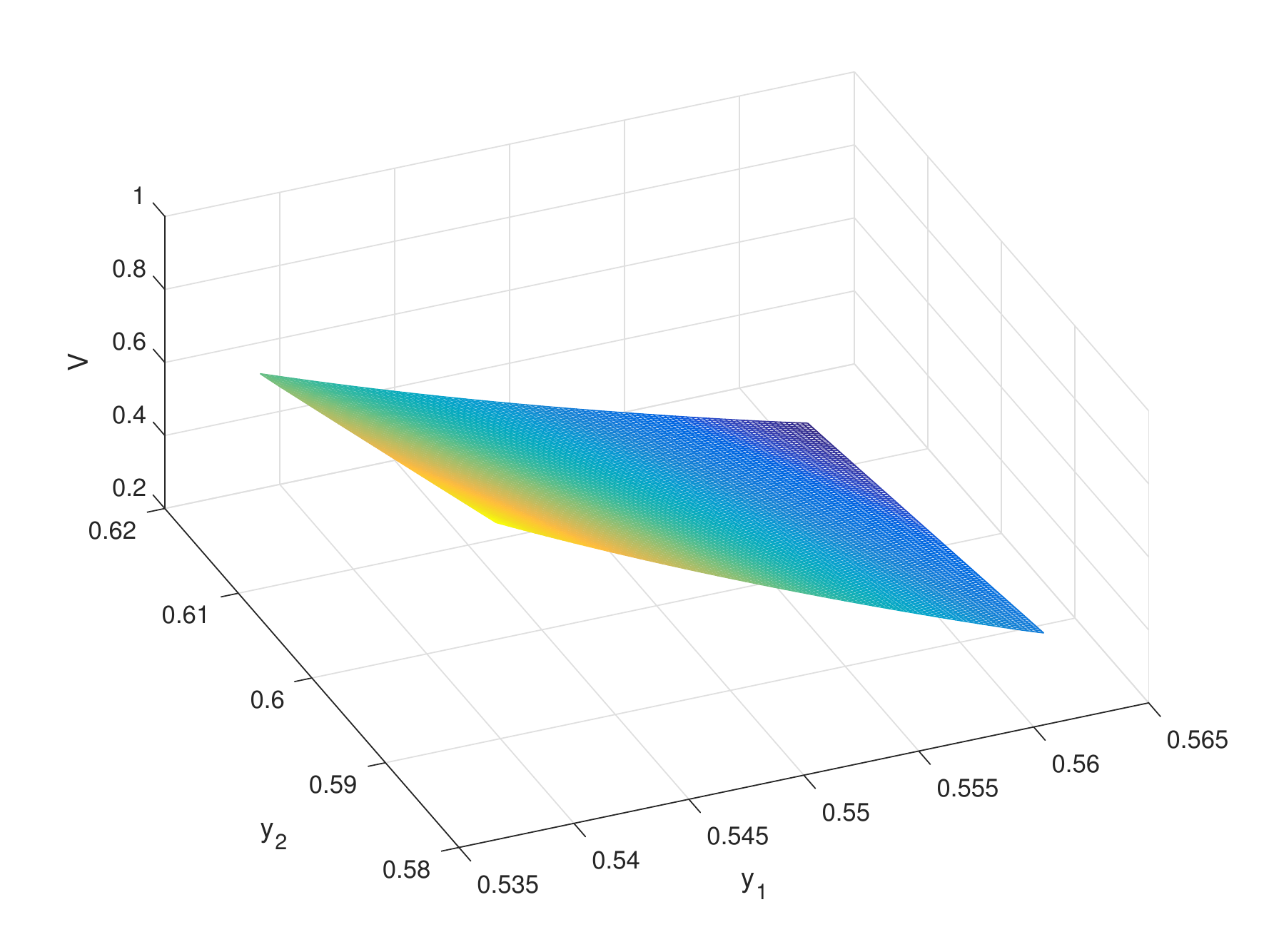}
  \caption{{\bf{Fig.10}} Normalised numerical solution $V$ at time step $10$ and $20$ respectively for $h=10$}
  \label{fig:3511}
\end{figure}  
\begin{figure}[h]
  \centering
  \includegraphics[ width=.43\textwidth]{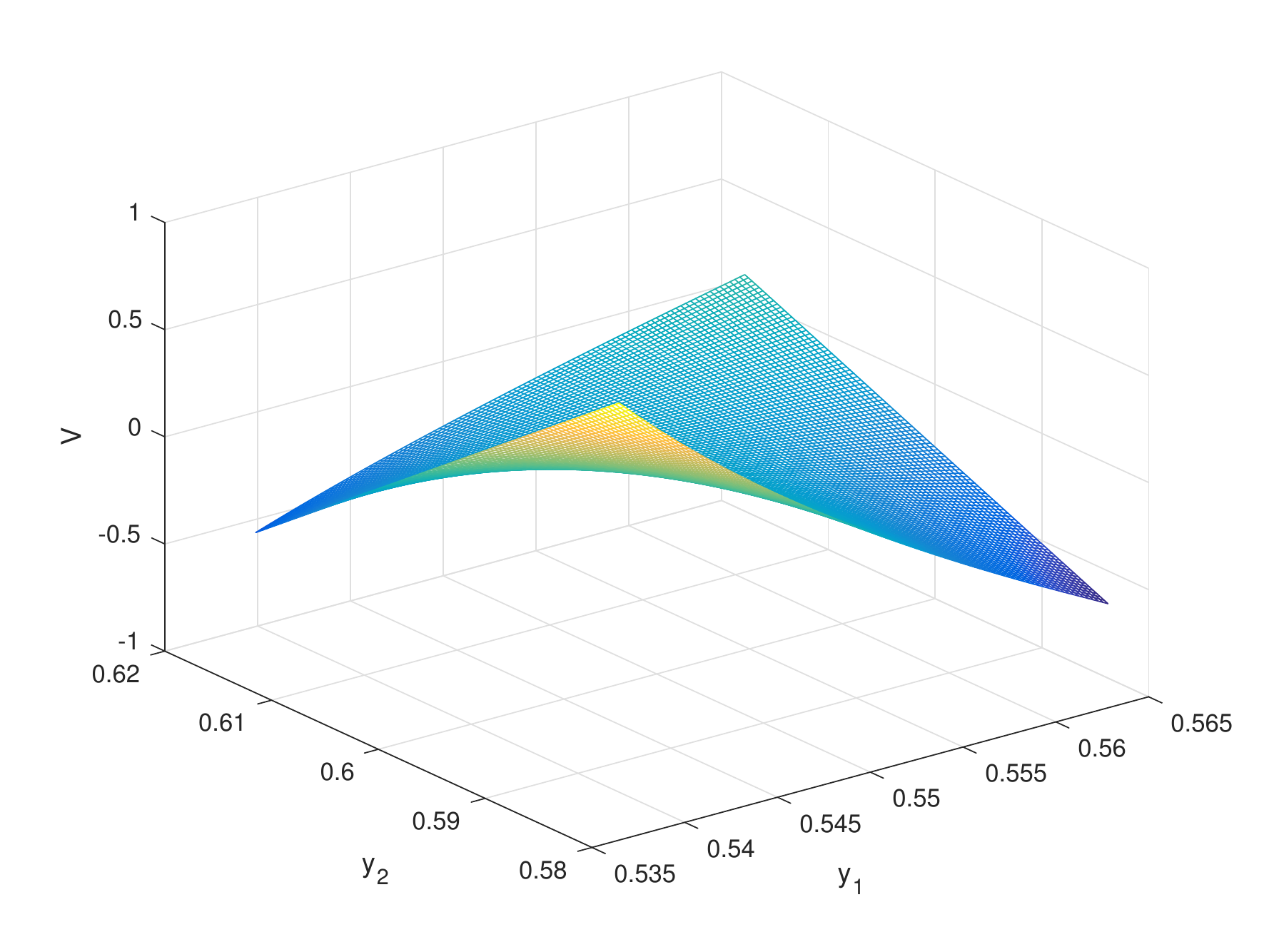}
  \includegraphics[ width=.43\textwidth]{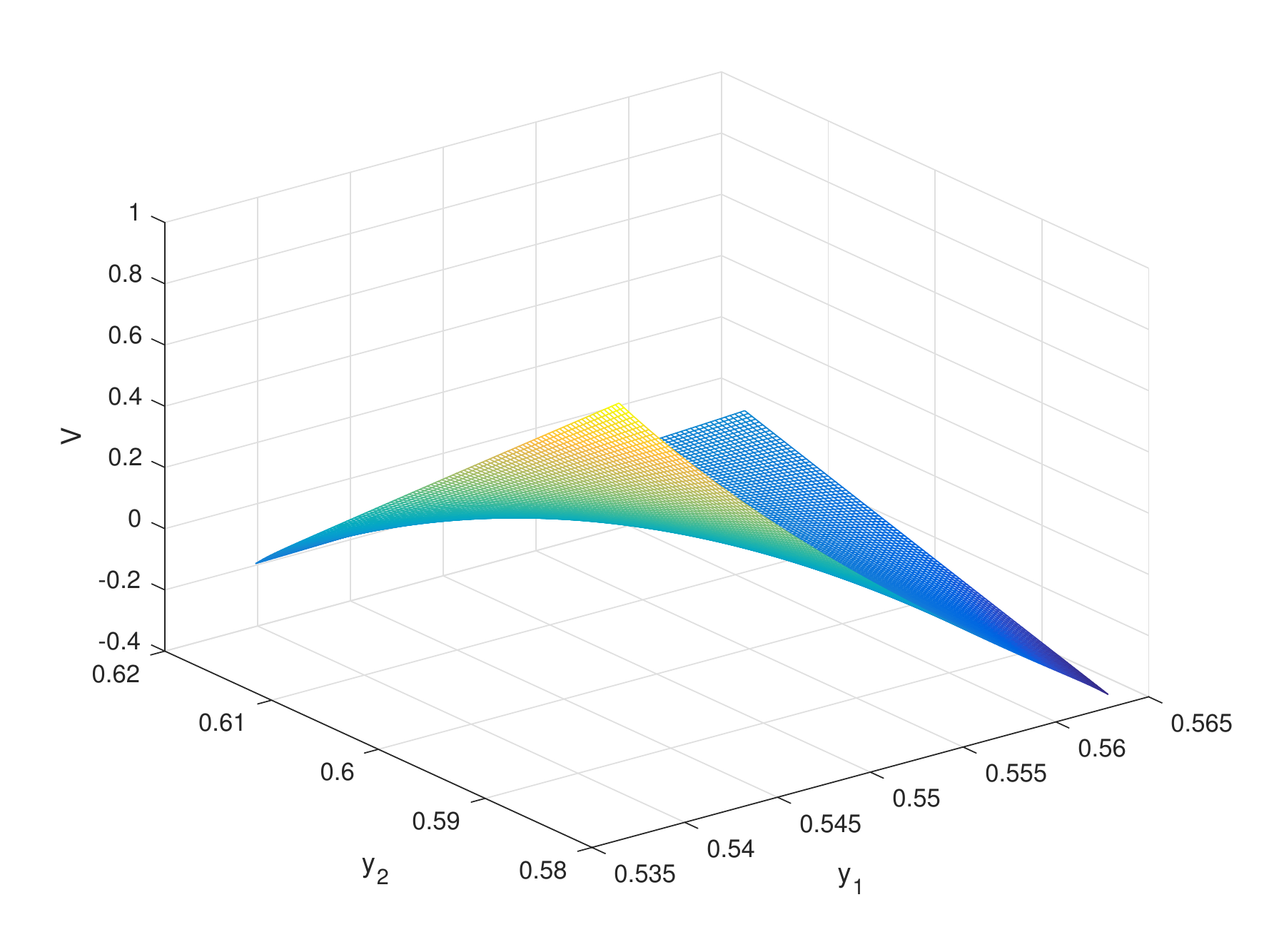}
  \caption{{\bf{Fig.11}} Normalised numerical  solution $V$ at time step $10$ and $20$ respectively for $h=50$}
  \label{fig:352}
\end{figure} 

\end{document}